\documentclass[ejsv2]{imsart}

\RequirePackage[numbers]{natbib}
\RequirePackage[colorlinks,citecolor=blue,urlcolor=blue]{hyperref}
\RequirePackage{graphicx}

\startlocaldefs
\theoremstyle{plain}

\newtheorem{theorem}{Theorem}[section]
\newtheorem{proposition}{Proposition}[section]
\newtheorem{remark}{Remark}[section]
\newtheorem{lemma}[theorem]{Lemma}

\theoremstyle{definition}

\theoremstyle{remark}

\endlocaldefs

\usepackage{lmodern}        
\usepackage{xcolor}
\usepackage{enumitem}
\usepackage{caption}
\usepackage{subcaption}
\usepackage{changepage}
\usepackage{comment}
\usepackage{algorithm}
\usepackage{algpseudocode}

\numberwithin{equation}{section}

\makeatletter
\@for\@tempa:={A,B,C,D,E,F,G,H,I,J,K,L,M,N,O,P,Q,R,S,T,U,V,W,X,Y,Z}\do{%
  \expandafter\edef\csname b\@tempa\endcsname{\noexpand\mathbb{\@tempa}}%
  \expandafter\edef\csname c\@tempa\endcsname{\noexpand\mathcal{\@tempa}}%
}
\makeatother

\begin{document}
\begin{frontmatter}
\title{Rate-optimal and computationally efficient nonparametric density and probability estimation on the circle and the sphere}
\runtitle{Nonparametric estimation on the circle and the sphere}

\begin{aug}
\author[A]{\fnms{Athanasios G.}~\snm{Georgiadis}\ead[label=e1]{georgiaa@tcd.ie}} 
\and
\author[A]{\fnms{Andrew P.}~\snm{Percival}\ead[label=e2]{aperciva@tcd.ie}} 

\address[A]{School of Computer Science and Statistics,
Trinity College of Dublin\printead[presep={,\ }]{e1,e2}}

\runauthor{A. G. Georgiadis and A. P. Percival}
\end{aug}

\begin{abstract}
We investigate the problem of density estimation on the unit circle and the unit sphere from a computational perspective. Our primary goal is to develop new density estimators that are both rate-optimal and computationally efficient for direct implementation. After establishing these estimators, we derive closed-form expressions for probability estimates over regions of the circle and the sphere. Then, the proposed theories are supported by extensive simulation studies. The considered settings naturally arise when analyzing phenomena on the Earth’s surface or in the sky (sphere), as well as directional or periodic phenomena (circle). The proposed approaches are broadly applicable, and we illustrate their usefulness through case studies in zoology, climatology, geophysics, and astronomy, which may be of independent interest. The methodologies developed here can be readily applied across a wide range of scientific domains.
\end{abstract}



\end{frontmatter}

\maketitle

\section{Introduction} 
Consider a random phenomenon expressed through some random variable $X$. The probability density function (PDF) $f=f_X$, or simply the \textit{density}, of $X$ compresses all vital information for the random variable $X$, and consequently the target phenomenon. Expectations, event probabilities, and related quantitative insights are fully determined by the knowledge of $f$.

Of course, the density is unknown, and the most we can do is to ``access" it from independent observations of it. In our context, we assume that $X_1,\dots,X_n$, $n\in\bN$ are independent realizations of $X$, which we must infer to assess the necessary information for $f$.

Nonparametric Statistics introduces \textit{estimators} of $f$, defined as functions of the form $\widehat{f}_n : \cM^{n} \times \cM \to \bR$, where $\cM$ denotes the set of all possible values that the random variable $X$ can take. Given the random sample $X_1, \dots, X_n$, we can construct 
$\widehat{f}_n(x) = \widehat{f}_n(x;X_1, \dots, X_n)$, $x \in \cM$,
which is referred to as the \textit{density estimator} of $f$. 

Nonparametric Statistics aims to provide such estimators without imposing restrictive---and often unrealistic---assumptions on the specific form of the unknown density. Instead, it seeks \textit{data-driven} procedures that consistently estimate $f$ under mild conditions, such as continuity or smoothness/regularity. Typical constructions rely on kernels, wavelets, splines, and other sophisticated mathematical tools. 

Historically, the earliest contributions are thanks to Rosenblatt~\cite{Ros}, Parzen \cite{Parzen}, and Breiman and Heller~\cite{BH}. Foundational monographs include \cite{HGPT,silverman}, while the book \cite{Tsybakov} is now regarded as a standard modern reference. For a representative, but in no sense exhaustive, list of further developments in the field, see~\cite{BT24,GNO1,GNO2,DG,DL2,DJKP,EY,GL,GL3,GLne1,GLne2,HWC,HI,IK,Lacroix-Lambert,KePT,KeLP,KLP2,Pel1,Pel2,RigT}, and the references therein .

In the present study we are interested in periodic random events and in random phenomena distributed on the Earth and Sky. The first can be modelled in terms of random variables distributed on the \textit{unit circle} $\bS^1=\{x\in\bR^2:\|x\|=1\}$, while the second are expressed via random variables distributed on the \textit{unit sphere} $\bS^2=\{x\in\bR^3:\|x\|=1\}$. Kernel density estimators (KDE) in that setting have a long history and applications in astrophysics; see for example \cite{BKMP}, \cite{BKMP_Asymptotics} and \cite{Geller2008}.

Recently, rate-optimal kernel density estimators in general metric spaces, covering the case of the spheres $\bS^d$, $d\in\bN$ have been established  and explored in \cite{CGKPP,CGP,CGW}.

The main challenge here is that the theoretically optimal estimators, established in the aforementioned works, are expressed in terms of infinite series. For example for the case of the sphere $\bS^2$, a ``smooth" density $f$, can be estimated by the KDE
\begin{equation}\label{estimator:second}
\widehat f_{n,h}(x):=\frac{1}{n}\sum_{j=1}^n\sum_{\ell=0}^\infty \frac{2\ell+1}{4\pi}g\Big(h\sqrt{\ell(\ell+1)}\Big)P_\ell(\langle x,X_j\rangle),\quad x\in\bS^2,
\end{equation}
for some appropriate function $g:\bR\to\bR$, and some bandwidth $h>0$. Note that $P_\ell$ in (\ref{estimator:second}) are the well-known Legendre polynomials; for more see Section \ref{sec: known estimators}. This automatically means that the KDE cannot be used by any software in its stand form. Therefore, when we are in front of a data-study, we need to address the following: 

\textit{Question:} Can we propose new estimators, that are computationally-efficient and remain rate-optimal?

The answer to the above question is \textit{affirmative}. The first problem we address in this study is to prove that there exists an index $N$ such that the finite-order estimator, obtained by summing from $\ell = 0$ to $\ell = N$, retains rate optimality. Note that this corresponds to a frequency cutoff, which is a natural practise in Approximation Theory. The corresponding results are provided in Section~\ref{sec: new estimators}.

\vspace{0.2cm}

Afterwards, with the established finite-order estimators in hand, we then turn our attention to \textit{probability estimation}. 

The \textit{new challenge} here is that the numerical evaluation of the integral of a given density estimator is relatively slow and computationally intensive. Therefore, in order to make these probability estimators more suitable for use, we analytically integrate these estimators over a given region into a closed expression. Consequently, numerical computation is required only for evaluating these closed forms. Section \ref{sec: probability} contains the derivation of the closed forms of these integrals and Figure \ref{fig: s1 s2 integration times} illustrates, at a glance, the resulting computational advantages.

\vspace{0.2cm}

In the sequel, we performed simulation studies to accompany these results. This can be found in Section \ref{sec: simulations} and showcases the accuracy of both the KDEs and probability estimators when compared to classic and mixture densities over the unit circle and sphere. This section also features evidence of the drastic time reduction in the evaluation of the closed forms of the probability estimators rather than the numerical integration method; see Section \ref{sec: comparison}.

Section~\ref{sec: case studies} presents four case studies in which our density estimators are applied to real-world data, with the aim not only of demonstrating the applicability of our theoretical results but also of disseminating the proposed methodology to the broader community. 

Our first case study, drawn from \textit{zoology}, concerns the movement directions of honeybees observed within a particular hive. As the data are naturally expressed as angles, our methods for the unit circle $\bS^1$ are directly applicable.

Our second case study, arising from \textit{climatology}, examines the annual precipitation patterns in Los Angeles, USA---a temporally periodic dataset. In such data, time can be interpreted as a periodic variable wrapped around the unit circle, making our methods on $\bS^1$ particularly suitable. More generally, any periodic phenomenon can be treated by our methodologies presented here.

Our third case study, originating from \textit{geophysics}, focuses on the spatial distribution of earthquake epicenters on the Earth’s surface. These locations are naturally represented as points on the unit sphere $\bS^2$, making the use of our spherical estimators particularly appropriate for analyzing global seismic activity.

Our fourth case study, drawn from \textit{astrophysics}, examines the distribution of bright stars in our galaxy. Each star is characterized by its galactic latitude and galactic longitude, providing data inherently defined on the unit sphere and thus ideally suited to our proposed estimators.

Finally, Section~\ref{sec: appendix} (the Appendix) contains the proofs of all new theoretical results established in this study.

\subsection*{Notation.} Through all the document the sets $\bN$, $\bN_0$, denote the positive and non-negative integers while $\bR$ and $\bC$, the real and complex numbers respectively. By $g\in\cC^{\tau}(\bR)$, we will mean that $g$ is a function in $\bR$ with continuous derivatives up to order $\tau\in\bN$. For any $s>0$, the number $\lfloor s\rfloor$ is the largest integer strictly smaller than $s$ and $\lceil s\rceil$ is the smallest integer strictly greater than $s$. We denote further $\cL^p(\bS^d)=\{f:\bS^d\to\bC:\|f\|_p<\infty\}$, $p\in[1,\infty]$, where $\|f\|_p^p=\int_{\bS^d}|f(x)|^p d\sigma_{\bS^d}(x)$, for $p<\infty$ and $\|f\|_{\infty}=\text{essup}|f(x)|$, the Lebesgue spaces on $\bS^d$. Finally $\langle x,y\rangle$ stands for the dot product of $x,y\in\bR^m$, $m\in\bN$.

\section{Background}
In this section, we collect the background material necessary for our study.

\subsection{The underlying geometries}\label{sec: geometries}
Firstly, let us formally introduce the settings of interest and their properties. 

We start with $\bS^1=\{x=(x_1,x_2)\in\bR^2:\|x\|=1\}$ the unit circle. Recall the polar representation of each $x=(x_1,x_2)=(\cos\theta,\sin\theta)=:x(\theta)$, where $\theta$ is the unique angle in $(-\pi,\pi]$ given by 
\begin{equation*}
\theta:=\begin{cases}
\arccos x_1,\quad x_2\ge0\\
-\arccos x_1,\quad x_2<0
\end{cases}.
\end{equation*}

The circle $\bS^1$ is naturally suitable for the study of periodic functions. The circular distance, is the angle between the corresponding points and the origin;
\begin{equation}
\nonumber
\rho_{\bS^1}(x,y)=\arccos(\langle x,y\rangle)=\min\{|\theta-\phi|,2\pi-|\theta-\phi|\},\; x=x(\theta),y=y(\phi)\in\bS^1.
\end{equation}
We denote by $d\sigma_{\bS^1}(x)$ the length measure in $\bS^1$. Then a measurable function $f:\bS^1\to\bC$ integrates as follows
\begin{equation*}
\int_{\bS^1}f(x)d\sigma_{\bS^1}(x)=\int_{-\pi}^{\pi}f(\cos\theta,\sin\theta) d\theta.
\end{equation*}
Finally, the geometry is equipped with the (negative of the) Laplacian, which is given by $L_{\bS^1}=-\frac{d^2}{d\theta^2}$.

We proceed to the sphere $\bS^2=\{x=(x_1,x_2,x_3)\in\bR^3:\|x\|=1\}$. We present the necessary notions released from \cite[Section 6]{Williams99} and \cite[Section 1]{DaiXu}. The \textit{ spherical coordinates} are given by
\begin{equation}
\label{spherical coordinates}
x_1=\sin\theta\cos\phi,\quad x_2=\sin\theta\sin\phi,\quad x_3=\cos\theta,\quad(\theta,\phi)\in[0,\pi]\times(-\pi,\pi],
\nonumber\end{equation}
and can be resolved as
\begin{equation}
\theta=\arccos(x_3),\quad\text{and}\quad \phi=\begin{cases}
    \pi,\quad&\text{if }\;x_1=x_2=0,\\
    \text{sign}(x_2)\arccos\Big(\frac{x_1}{\sqrt{x_1^2+x_2^2}}\Big),\quad &\text{otherwise}
\end{cases}.
\nonumber\end{equation}
With a slight abuse of the notation we sometimes denote a point $x\in\bS^2$ just by its spherical coordinates $x=(\theta,\phi)$.

We denote by $d\sigma_{\bS^2}(x)$ the surface measure in $\bS^2$. Then a measurable function $f:\bS^2\to\bC$ integrates as follows
\begin{equation*}
\int_{\bS^2}f(x)d\sigma_{\bS^2}(x)=\int_{-\pi}^{\pi}\int_{0}^{\pi}f(\sin\theta\cos\phi,\sin\theta\sin\phi,\cos\theta)\sin\theta d\theta d\phi.
\end{equation*}
The spherical distance is defined, as in $\bS^1$, by
\begin{equation}\label{metric on S2}
\nonumber
\rho_{\bS^2}(x,y)=\arccos(\langle x,y\rangle),\quad x,y\in\bS^2.
\end{equation}
The (negative of the) spherical Laplacian in spherical coordinates is given by
\begin{equation*}
L_{\bS^2}=-\frac{1}{\sin\theta}\frac{\partial}{\partial\theta}\Big(\sin\theta\frac{\partial}{\partial\theta}\Big)-\frac{1}{\sin^2\theta}\frac{\partial^2}{\partial^2\phi}.
\end{equation*}

\subsection{Legendre polynomials and spherical harmonics}\label{sec: spherical} Among the most valuable tools on the sphere, are the spherical harmonics and Legendre polynomials. 

The \textit{Legendre polynomials} $P_{\ell}$, $\ell\in\bN_0$, are given by the simple Rodrigues formula
\begin{equation}
\label{Rodrigues}
P_{\ell}(u)=\frac{1}{2^{\ell}\ell!}\frac{d^{\ell}}{du^{\ell}}(u^2-1)^{\ell},\quad u\in[-1,1],\ell\in\bN_0.
\end{equation}
We will further need the next precise expression of the Legendre polynomials:
\begin{equation}
\label{Simple for Legendre polynomials}
P_{\ell}(u)=\sum_{k=0}^{\ell}\binom{\ell}{k}\binom{\ell+k}{k}\frac{(-1)^k}{2^k}(1-u)^k,\quad u\in[-1,1],\ell\in\bN_0,
\end{equation}
and the well-known upper bound
\begin{equation}
\label{legendre maximum}
|P_{\ell}(u)|\le 1,\quad u\in[-1,1],\ell\in\bN_0.
\end{equation}
Note further, that $P_{\ell}(1)=1$, for every $\ell\in\bN_0$, $P_{0}\equiv1$ and the set $\{P_{\ell}:\ell\in\bN_0\}$ consists of an orthogonal basis for $L^2[-1,1]$, the space of square integrable functions on $[-1,1]$.

The \textit{associated Legendre functions} $P^m_{\ell}:\ell\in\bN_0,-\ell\le m\le m$, are defined by

\begin{equation}\label{associated Legendre functions}
P_{\ell}^{m}(u):=\begin{cases}
(-1)^m(1-u^2)^{m/2}\frac{d^m}{du^m}P_{\ell}(u),\quad m\ge0
\\
(-1)^m\frac{(\ell+m)!}{(\ell-m)!}P_{\ell}^{-m}(u),\quad m<0
\end{cases},\quad \ell\in\bN_0.
\end{equation}
Of course $P_{\ell}^{0}=P_{\ell}$, for every $\ell\in\bN_0$.

Finally, for every $\ell\in\bN_0$ and $0\le m\le \ell$ the spherical harmonic $Y^{m}_{\ell}:[0,\pi]\times(-\pi,\pi]\to\bC$ are given by
\begin{align}
\label{spherical by associated Legendre}
Y_{\ell}^m(\theta,\phi)&=\sqrt{\frac{2\ell+1}{4\pi}}\sqrt{\frac{(\ell-m)!}{(\ell+m)!}}P_{\ell}^{m}(\cos\theta)e^{im\phi}
=:N_{\ell,m}P_{\ell}^{m}(\cos\theta)e^{im\phi},
\end{align}
and are expanded for $-\ell\le m<0$ by the following conjugation; 
\begin{equation}
\label{spherical harmonics conjugate}
Y^{m}_{\ell}(\theta,\phi)=(-1)^m\overline{Y^{-m}_{\ell}(\theta,\phi)},\quad \ell\in\bN_0,\;-\ell\le m\le 0.
\end{equation}
The most important property of spherical harmonics is that their union $\{Y^m_{\ell}:\ell\in\bN_0,-\ell\le m\le \ell\}$ consists of an orthonormal basis of $(\cL^2(\bS^2),d\sigma_{\bS^2})$. Moreover, they are connected with Legendre polynomials under the fundamental additive formula:
\begin{equation}
\label{Legendreand spherical harmonics}
\frac{2\ell+1}{4\pi}P_\ell(\langle x,y\rangle)=\sum_{m=-\ell}^{\ell}Y^m_{\ell}(\theta_x,\phi_x)\overline{Y_{\ell}^m(\theta_y,\phi_y)},\quad x,y\in\bS^2,\;\ell\in\bN_0.
\end{equation}

\subsection{Density estimation on the circle and sphere}\label{sec: known estimators}
Let us now summarise some key results from \cite{CGKPP},\cite{CGP}, and \cite{CGW}, in which density estimation has been studied on metric spaces $\cM$ associated with operators $L$ under mild assumptions. In our case, the spaces $\bS^d$, for $d=1,2$, with corresponding operators, the negative Laplacians $L_{\bS^d}$, given in Section \ref{sec: geometries}, obey these assumptions and therefore the results from those studies hold automatically. We present those results in the case of $\bS^d$, for $d=1,2$.

Let $g:\bR\to\bR$ be a measurable function such that $g(0)=1$ with a ``fast" decay (this is specified in Theorem \ref{th: milk metrika}). Such a function is referred to as a \textit{symbol} and will generate a KDE, as in the case of kernel functions on $\bR^d$. Let $h>0$ be a quantity that extends the notion of a bandwidth. 

Consider a random variable $X$ distributed on $\bS^d$, with unknown density $f$ and $X_1,\dots,X_n$, $n\in\bN$ iid copies of $X$. 

Let $d=1$. The KDE $\widehat f_{n,h}$ takes the following form:
\begin{equation}
\label{eq: kde on S1}
    \widehat f_{n,h}(x)=\frac{1}{2\pi n}\sum_{j=1}^n\Big(1+2\sum_{\ell=1}^\infty g(h \ell)\cos(\ell(\theta-\theta_j))\Big),\quad x=x(\theta),\;\theta\in(-\pi,\pi],
\end{equation}
where $\theta_j$ the polar angles of $X_j$, $j=1,\dots,n$ and $\theta$ the polar angle of the arbitrary point $x\in\bS^1$.

Let $d=2$. The KDE $\widehat f_{n,h}$ takes the following form:
\begin{equation}
\label{eq: kde S2}
    \widehat f_{n,h}
(x)=\frac{1}{n}\sum_{j=1}^n\sum_{\ell=0}^\infty\frac{2\ell+1}{4\pi}g\Big(h\sqrt{(\ell(\ell+1))}\Big)P_\ell(\langle x,X_j\rangle),\quad x\in\bS^2.
\end{equation}

As always in Nonparametric statistics, we assume that the unknown density, lies on a regularity space. Let $s>0$. We say that a function $f$ belongs to the H\"{o}lder space of regularity order $s$, $\dot{\cH}^s(\bS^d)$, when 
\begin{equation}
\label{Holder norm}
\|f\|_{\dot{\cH}^s(\bS^d)}=\|f\|_{\dot{\cH}^s}=\sup_{x\neq y}\frac{|L_{\bS^d}^{k/2}f(x)-L_{\bS^d}^{k/2}(y)|}{\rho_{\bS^d}(x,y)^{s-k}}<\infty,
\nonumber\end{equation}
where $k$ the smallest integer strictly less than $s$.

Observe that $L_{\bS^d}$ are second order differential operators, which indicates their use for defining regularity spaces and the half in the exponent. For the fractional powers of the operators, we refer the reader to \cite[p. 263]{RS}. For more details on regularity spaces on the sphere and more general spaces, we refer to \cite{CKP,GKembeddings,GN,KOPP,KP}

Finally, the simplest measure of the accuracy of an estimator is the mean squared error:
\begin{align}\nonumber
    \text{MSE}(&\widehat f_{n,h})(x):=\bE_f[(\widehat f_{n,h}(x)-f(x))^2]\\ &=\int_{\bS^d} \dots \int_{\bS^d} (\widehat f_{n,h}(x;x_1,\dots,x_n)-f(x))^2 f(x_1)\dots f(x_n)\times \nonumber\\
    &\hspace{3cm}\times d\sigma_{\bS^d}(x_1)\dots d\sigma_{\bS^d}(x_n).
\end{align}

We now fix some terminology that will be used throughout the paper. 
Let $s>0$. We say that a function $g:\mathbb{R}\to\mathbb{R}$ is an \emph{$s$-admissible symbol} if it is even and measurable, belongs to $\mathcal{C}^{\tau}(\mathbb{R})$ for some $\tau>d+s$, satisfies $g(0)=1$ and $g^{(\nu)}(0)=0$ for all $1\le \nu\le \tau$, and if there exist constants $r>\tau+d$ and $C_\tau>0$ such that
\begin{equation}\label{eq: decay}
|g^{(\nu)}(\lambda)| \le C_\tau (1+|\lambda|)^{-r},
\qquad \text{for all } \lambda\in\mathbb{R}, \; 0\le \nu\le \tau.
\end{equation}
By Theorem 3.1 in \cite{CGW} we infer:

\begin{theorem}\label{th: milk metrika}
Consider a random variable $X$ distributed on $\bS^d$, $d\in\{1,2\}$, with unknown density $f$ and $X_1,\dots,X_n$, $n\in\bN$ iid copies of $X$. Assume that $f\in\cL^\infty(\bS^d)\cap\dot\cH^s(\bS^d)$ for some $s>0$. Let also $h=h_n=n^{-1/(2s+d)}$ and $g$ and $s$-admissible symbol.

Then the corresponding KDEs $\widehat f_{n,h}$ as in \eqref{eq: kde on S1}, for $d=1$, and 
\eqref{eq: kde S2} for $d=2$ satisfy,
\begin{equation}
\sup_{x\in\bS^d} \textnormal{MSE} (\widehat f_{n,h})(x)\leq cC(f) n^{-2s/(2s+d)},
\end{equation}
where $c>0$ depends only on $\tau,s,C_\tau$ and $d$, while $C(f)=\max\{\|f\|_\infty,\|f\|^2_{\dot \cH^s}\}$.
\end{theorem}

The reader is referred to \cite{CGKPP,CGP} for more general results in terms of errors and function spaces and to \cite[Theorem 11]{BKMP} for the lower bounds and thus the optimality of the rate above. 

The assumptions on the symbol are those required to generalize the usual notion of ``kernels of order $s$'', as given, for example, in \cite[Definition 1.3]{Tsybakov}, to general metric spaces. The details can be found in \cite[Section 3]{CGW}.

\section{Rate-optimal and computationally efficient estimators}\label{sec: new estimators} 

Whilst Theorem \ref{th: milk metrika} shows that the KDEs in \eqref{eq: kde on S1} and \eqref{eq: kde S2} are theoretically suitable estimators, the implementation of these KDEs comes with a computational challenge; each expression contains an infinite sum over $\ell$. To evaluate these estimators in practise, we must replace the infinite series by a finite-order version, keeping only terms up to some cutoff index $N$. Determining this index $N$ in such a way that the estimator \textit{remains rate-optimal} motivates us to establish the following theorem, whose proof, as all proofs of our study, can be found in the Appendix:

\begin{theorem}\label{th: finite sum}
Consider a random variable $X$ distributed on $\bS^d$, $d\in\{1,2\}$, with unknown density $f$ and $X_1,\dots,X_n$, $n\in\bN$ iid copies of $X$. Assume that $f\in\cL^\infty(\bS^d)\cap\dot\cH^s(\bS^d)$ for some $s>0$. Let also $h=h_n=n^{-1/(2s+d)}$ and $g$ and $s$-admissible symbol. We denote by 
\begin{equation}
\label{eq: terminal N Sd}
N_s := \Big\lfloor
c_{rd}\, n^{\frac{s+r}{(2s+d)(r-d)}}
\Big\rfloor + 1,
\end{equation}
where $c_{rd}=(d\pi(r-d))^{-\frac{1}{r-d}}$.

We introduce the following ---finite-order--- estimators:

$\bullet$ When $d=1$ 
\begin{equation}
\label{eq: kde S1 finite}
\widehat{f}_s(x):=\frac{1}{2\pi n}\sum_{j=1}^n\Big(1+2\sum_{\ell=1}^{N_s} g(h\ell)\cos(\ell(\theta-\theta_j))\Big),\quad x=x(\theta)\in\bS^1,
\end{equation}
where $\theta,\theta_j\in(-\pi,\pi]$ represent the polar angles of the arbitrary $x\in\bS^1$ and $X_j$, $j=1,\dots,n$, respectively.

$\bullet$ When $d=2$ 
\begin{equation}
\label{eq: kde S2 finite}
\widehat f_s(x):=\frac{1}{n}\sum_{j=1}^n\sum_{\ell=0}^{N_s} \frac{2\ell+1}{4\pi}g\big(h\sqrt{\ell(\ell+1)}\big)P_\ell(\langle x, X_j\rangle),\quad x\in\bS^2.
\end{equation}
Then the estimator $\widehat f_s$ as in (\ref{eq: kde S1 finite}) for $d=1$ and in (\ref{eq: kde S2 finite}) for $d=2$ enjoys
\begin{equation}
\label{eq: s2 trunc thm bound}
\sup_{x\in\bS^d}\textnormal{MSE}(\widehat{f}_s)(x)\leq C' n^{-2s/(2s+d)},
\end{equation}
where $C'$ depends on $C_{\tau}$ as well as $c$ and $C(f)$ as in Theorem \ref{th: milk metrika}.
\end{theorem}

We close the present section by listing the following remarks.

\begin{remark}\label{Remark 1}

(a) Our finite-order KDES in (\ref{eq: kde S1 finite}) and (\ref{eq: kde S2 finite}) of Theorem \ref{th: finite sum}, are computationally efficient, once we fix an appropriate $s$-admissible symbol $g$. Let $d\in\{1,2\}$ and $s>0$ and set $r:=2d+\lceil s\rceil+1$. Then the function
\begin{equation}
\label{gr symbol}
g_{r}(\lambda):=\frac{1}{1+|\lambda|^r},\quad\lambda\in\bR,
\end{equation}
is an $s$-admissible symbol \cite{CGP,CGW}. This means that our new optimal estimators are ready to be used in data analysis. This completes the line of work developed in \cite{CGKPP,CGP,CGW}, now complemented by the computational and applied aspects presented in the current study. For simulations and applications in the current study we refer the reader to Sections \ref{sec: simulations} and \ref{sec: case studies}.

(b) The finite estimators in Theorem \ref{th: finite sum}, indeed integrate to the unit; $\int_{\bS^d}\widehat{f}_s(x)d\sigma_{\bS^d}(x)=1$, as one would expect. For the case of the circle $\bS^1$, we can confirm this just by integrating (\ref{eq: kde S1 finite}) and using the properties of the trigonometric functions. For the case of the sphere $\bS^2$, we just have to use the Funk–Hecke formula; \cite[Theorem 1.2.9]{DaiXu} to obtain for every $y\in\bS^2$,
\begin{equation*}
\int_{\bS^2}P_{\ell}(\langle x,y\rangle)d\sigma_{\bS^2}(x)=2\pi\int_{-1}^{1}P_{\ell}(u)du=\begin{cases}
4\pi,\quad\ell=0
\\
0, \quad \ell\neq0
\end{cases},
\end{equation*}
with the last holding because $P_{0}\equiv1$ and thanks to the orthogonality of the Legendre polynomials.

(c) The KDEs obtained here, in the language of Harmonic Analysis, are commonly referred as band-limited. This is because they are obtained by Fourier expansions, using a finite part of the spectrum of the associated operators (Laplacians). This is a very standard direction in Analysis and Approximation Theory and finds applications in a wide range of areas.

(d) For clarity, the stopping index $N_s$, defined in (\ref{eq: terminal N Sd}), depends on both the sample size $n$ and the regularity parameter $s$.
\end{remark}

\section{Probability estimation}\label{sec: probability}
A problem of fundamental interest is the estimation of the probability of an event occurring within a certain region of a space. Our kernel density estimators can be used for such an application, not only as plug-in estimators but also as a rate optimal estimation for such a probability on the circle or sphere. We introduce the following probability estimators for every measurable set $A\subset\bS^d$, $d=1,2$ and we prove that these can indeed consistently estimate the probability of any measurable subset $A$.

\begin{proposition}\label{Th:MSEboundProbabilityEstimator}
Consider a random variable $X$ distributed on $\bS^d$, $d\in\{1,2\}$, with unknown density $f$ and $X_1,\dots,X_n$, $n\in\bN$ iid samples of $X$. Assume that $f\in\cL^\infty(\bS^d)\cap\dot\cH^s(\bS^d)$ for some $s>0$ and denote by $\widehat{f}_s$ the finite KDE as in Theorem \ref{th: finite sum}. Let finally $A\subset\bS^d$ be a measurable set. Then the following probability estimator 
\begin{equation}\label{eq: probabilityestimator}
\widehat{\bP}(A):=\widehat{\bP}(A,X_1,\dots,X_n):=\int_{A}\widehat{f}_s(x) d\sigma_{\bS^d}(x),
\end{equation}
satisfies
\begin{equation*}
\bE_f[(\widehat{\bP}(A)-\bP(A))^2]\leq c'n^{-2s/(2s+d)},
\end{equation*}
where $c'>0$ depends on the measure of $\sigma_{\bS^d}(A)$ as well as $C'$ from Theorem \ref{th: finite sum}.
\end{proposition}

With Proposition \ref{Th:MSEboundProbabilityEstimator} in hand we can integrate our new estimators and obtain probabilities of certain regions of $\bS^d$. We start with the case of $d=1$. Let $-\pi<\vartheta_1<\vartheta_2\le\pi$ and we introduce the following angular domain
\begin{equation}
\label{eq: region on S1}
A[\vartheta_1,\vartheta_2]=\{x=x(\theta)\in\bS^1:\;\vartheta_1\le\theta\le\vartheta_2\}.
\end{equation}

We apply Proposition \ref{Th:MSEboundProbabilityEstimator} to get the following:
\begin{proposition}
\label{th: integral S1 kde}
Let $X$ a random variable distributed on $\bS^1$,  with density $f\in\cL^\infty(\bS^1)\cap\dot\cH^s(\bS^1)$ for some $s>0$ and $X_1,\dots,X_n$ iid samples of $X$. Let also $h=h_n=n^{-1/(2s+d)}$ and $g$ and $s$-admissible symbol. Then
\begin{align}
\label{eq: kde S1 integrated}
\widehat{\bP}(A[\vartheta_1,\vartheta_2])&=\frac{\vartheta_2-\vartheta_1}{2\pi }
\\
&+\frac{1}{\pi n}\sum_{j=1}^n\sum_{\ell=1}^{N_s} \frac{g(h\ell)}{\ell}\left(\sin(\ell(\vartheta_2-\theta_j))-\sin(\ell(\vartheta_1-\theta_j))\right),
\nonumber\end{align}
where $\theta_j$ the polar angles of $X_j$, $j=1,\dots,n$.
\end{proposition}
\vspace{0.2cm}

The integration of the finite estimators (\ref{eq: kde S1 finite}) on the circle, is straightforward; see the Appendix. In contrast, our estimators (\ref{eq: kde S2 finite}) for $\bS^2$ involve the Legendre polynomials and their integration turns out to be challenging. Numerical methods end up taking a very long time to compute integrals of the form
$$\int_{A}P_\ell(\langle x, x_j\rangle)d\sigma_{\bS^2}(x),\quad A\subset\bS^2,$$
no matter how simple form the domain $A$ has, and therefore the software suffers providing the final estimated probabilities.

This motivates the direct computation of the integral of the finite estimator \eqref{eq: kde S2 finite}, as we will do in the next theorem. We will see that this closed form is also particularly suitable for vectorisation when being evaluated, further decreasing run-time, especially when combined with parallel programming. The anti-derivatives of the form of the Legendre polynomials in \eqref{eq: kde S2 finite} can be written in a closed form but special functions will be needed in this expression. These special functions are well-known and are found natively in software like \textit{Mathematica} or are easily available through standard packages, such as the \texttt{scipy} package for \textit{Python}.

We will particularly need the so-called \textit{incomplete Beta function}.  For every $a,b>0$ the incomplete Beta function is defined by
\begin{equation}
\label{incomplete B function}
B(u;a,b):=\int_{0}^{u}t^{a-1}(1-t)^{b-1}dt,\quad \text{for every }u\in[0,1].
\end{equation}
Note that for $u=1$, its value is equal to the standard Beta function. Furthermore, if $a,b\in\bN$, then the integrand function in (\ref{incomplete B function}) is a polynomial, so the incomplete Beta function, its anti-derivative, is just a polynomial too. For the better presentation of the study we introduce the following kernels: For every $\ell\in\bN$, $1\le m\le k\le\ell$ and $\gamma_1,\gamma_2\in[-1,1]$, we denote by
\begin{equation}\label{cBkm}
\cB_{m,k}(\gamma_1,\gamma_2):=B\Big(\frac{1+\gamma_2}{2};\frac{m}{2}+1,k-\frac{m}{2}+1\Big)-B\Big(\frac{1+\gamma_1}{2};\frac{m}{2}+1,k-\frac{m}{2}+1\Big).
\end{equation}

As in the case of the circle we work with angular domains; namely denote $\vartheta=(\vartheta_1,\vartheta_2),\;\varphi=(\varphi_1,\varphi_2)$, where $0\le\vartheta_1<\vartheta_2\le\pi$ and $-\pi\le\varphi_1<\varphi_2\le\pi$. We set
$$\cA_{\vartheta,\varphi}:=\{x=x(\theta,\phi)\in\bS^2:(\theta,\phi)\in[\vartheta_1,\vartheta_2]\times[\varphi_1,\varphi_2]\}.$$

\begin{theorem}\label{th: probability S2}
Let $X$ a random variable distributed on $\bS^2$,  with density $f\in\cL^\infty(\bS^2)\cap\dot\cH^s(\bS^2)$ for some $s>0$ and $X_1,\dots,X_n$ iid samples of it. Let also $h=h_n=n^{-1/(2s+d)}$ and $g$ and $s$-admissible symbol. Then for every $\vartheta=(\vartheta_1,\vartheta_2),\;\varphi=(\varphi_1,\varphi_2)$, with $0\le\vartheta_1<\vartheta_2\le\pi$ and $-\pi\le\varphi_1<\varphi_2\le\pi$, we have that
\begin{equation}
\label{eq: kde s2 integrated}
\widehat{\bP}(\cA_{\vartheta,\varphi})=\frac{1}{n}\sum_{j=1}^n\sum_{\ell=0}^{N_s} g\Big(h\sqrt{\ell(\ell+1)}\Big)\Big(I_{\ell}^{0}(j)+2\sum_{m=1}^{\ell}I_{\ell}^{m}(j)\Big),
\end{equation}
where:

for every $\ell\in\bN_0$ and $j=1,\dots,n$,
\begin{align*}
I_{\ell}^{0}(j) = N_{\ell,0}^2 (\varphi_2 - \varphi_1) P_{\ell}(\cos\theta_j) 
\sum_{k=0}^{\ell} 
& { \binom{\ell}{k} \binom{\ell + k}{k} }
\left(\frac{-1}{2} \right)^k\cdot \frac{1}{k+1}
\\
&\left( (1 - \cos\vartheta_2)^{k+1} - (1 - \cos\vartheta_1)^{k+1} \right);
\end{align*}

for every $\ell\in\bN$, $1\le m\le \ell$ and $j=1,\dots,n$,
\begin{align*}
I^{m}_\ell (j)&=\frac{1}{m}N_{\ell,m}^2P^m_{\ell}(\cos\theta_j)\Big(\sin\big(m(\varphi_2-\phi_j)\big)-\sin\big(m(\varphi_1-\phi_j)\big)\Big)
\times
\\
&\hspace{0.6cm}\times\sum_{k=m}^{\ell}{ \binom{\ell}{k}\binom{\ell+k}{k}}\frac{2(-1)^k k!}{(k-m)!}\cB_{k,m}(\cos\vartheta_2,\cos\vartheta_1)
\end{align*}
for $N_{\ell,m}$ as in (\ref{spherical by associated Legendre}), $\cB_{k,m}$ as in (\ref{cBkm}) and where $(\theta_j,\phi_j)$ are the spherical coordinates of the observations $X_j$, $j=1,\dots,n$.
\end{theorem}

The proof of the above theorem, relies on the next result, which is of potential independent interest, beyond the field of Statistics.

\begin{lemma}
\label{Lemma: integrals related to incomplete Beta function}
Let $-1\le \gamma_1<\gamma_2\le 1$, $\ell\in\bN$ and $1\le m\le\ell$. Then the associated Legendre functions $P_\ell^m$ integrate in the closed interval $[\gamma_1,\gamma_2]$ as
\begin{align}
\int_{\gamma_1}^{\gamma_2}P_{\ell}^m(u)du&=\sum_{k=m}^{\ell}\binom{\ell}{k}\binom{\ell+k}{k}\frac{2(-1)^k k!}{(k-m)!}
\cB_{k,m}(\gamma_1,\gamma_2)
\nonumber
\end{align}
where $\cB_{k,m}$ the kernels in (\ref{cBkm}) involving the incomplete Beta function (\ref{incomplete B function}).
\end{lemma}

\section{Simulations}\label{sec: simulations}

In this section, we reinforce with simulations the theoretical finds in Sections \ref{sec: new estimators} and \ref{sec: probability}. We simulate from classic distributions on the sphere and circle and present our estimation of the densities and probabilities; that is the implementation of Theorem \ref{th: finite sum}, Proposition \ref{th: integral S1 kde} and Theorem \ref{th: probability S2}. In particular, we first sample from the uniform distributions on $\bS^d$ and then the very well-known von Mises-Fisher distribution.

Our methodology for this section is as follows: we select a random variable on $\bS^1$ or $\bS^2$ with a density $f$ and simulate $n=1000$ points from its distribution. We then construct finite-order kernel density estimators $\widehat {f}_s$, as in Theorem \ref{th: finite sum}, taking $s=0.5$, $1$, and $2$ with $r=2d+\lceil s\rceil+1$ and symbol $g_r$, as in Remark \ref{Remark 1}. We then plot the original distribution $f$ alongside the simulated sample of $n=1000$ points and the KDEs $\widehat f_{0.5}$, $\widehat f_1$, and $\widehat f_2$. 

In the case of $\bS^1$, the plot is given on an annulus (with an inner radius of $0.8$ and outer radius of $1.2$) with the colour representing the value of the function at a given angle, achieved with the \texttt{ColorFunction} option in Mathematica's \texttt{PolarPlot}. The points are also given on the same annulus with the angle of each point in polar coordinates being equal to the angle of the observation but with the radius given as a random number between $0.8$ and $1.2$ in order to better spread out the points for visual clarity. 

In the case of $\bS^2$, the plot is given on a 3D rendering of a sphere, viewed from the point $(1.3,-2.4,2)$, the default in Mathematica's \texttt{SphericalPlot3D}. The spheres are coloured according to their corresponding function (the PDF or KDE). This colour plot was achieved by evaluating each function over a \texttt{Table} in terms of $\theta$ and $\phi$. In these plots we took the step size $\pi/32$ in each coordinate giving a table of $65\times33$ points which were entered into \texttt{ListDensityPlot} to get an \texttt{Image} which was used as the \texttt{Texture} in the \texttt{PlotStyle} option of \texttt{SphericalPlot3D}. The $n=1000$ points are given on the same sphere but with no added colour.

For each simulated KDE, we estimated the probability over several regions using the methods established in Proposition \ref{th: integral S1 kde} and Theorem \ref{th: probability S2} for $\bS^1$ and $\bS^2$, respectively. Precisely, the Tables \ref{tab:unif probs}-\ref{tab:vm mixed probs} contain: these estimated probabilities compared to the observed frequency in the given region as well as the true (theoretical) probability, obtained by integrating the known PDF over the same region. 

Finally, alongside the probabilities in these tables, we also include an indicator of the performance of each KDE: the estimated mean integrated square error (MISE). The MISE is a global criterion for the fit of an estimator $\widehat f$ to a density $f$ and is given by
\begin{equation}
    \text{MISE}(\widehat f)=\bE_f \Big[\int_{\bS^d} (\widehat f(x)-f(x))^2 d\sigma_{\bS^d}(x)\Big],\quad d=1,2.
\end{equation}
Note also that, by the Fubini-Tonelli theorem and Theorem \ref{th: finite sum}, we have when $f\in\cL^{\infty}(\bS^d)\cap\dot{\cH}(\bS^d)$, for some $s>0$ and $\widehat{f}_s$ as in Theorem \ref{th: finite sum},
\begin{equation}
    \text{MISE}(\widehat {f}_s)\leq \sigma_{\bS^d}(\bS^d) C' n^{-2s/(2s+d)},\quad d=1,2.
\end{equation}
In order to estimate the MISE, we simulated $n=1000$ points, $30$ times, computed the MISE in each instance and took the mean of this sample. 

The form of the results of Proposition \ref{th: integral S1 kde} and Theorem \ref{th: probability S2} are suitable for array programming. This means that for a given region of $\bS^1$ or $\bS^2$ to integrate over, one can apply functions such as the associated Legendre polynomials $P_m^\ell$ over the entire dataset for all necessary values of $m$ and $\ell$ simultaneously. This drastically reduces computation time, especially when compared to numerical integration (see Section \ref{sec: comparison}). Finally, a computational tip is that there is no need to calculate the inclination angles $\theta_j$ as this only appears as the argument of $\cos$ and by (\ref{spherical coordinates}) $\cos\theta_j=X_{j3}$, the third coordinate of the $j^{\text{th}}$ observation. 

We will study $6$ scenaria presented below.

\subsection{The uniform distribution}
\label{sec: unif dist}

\paragraph*{Scenario 1 (Uniform distribution on $\mathbb{S}^2$).}
We begin with simulations of the uniform distribution on the unit sphere. In order to simulate from this distribution on $\bS^d$, the usual approach, from \cite{Muller1956} is given in Algorithm \ref{alg: unif Sd}.
\begin{algorithm}
\caption{Simulation of uniformly distributed points on $\bS^d$}\label{alg: unif Sd}
\begin{algorithmic}[1]
    \Require Dimension of the sphere: $d\in\bN$
    \Require Sample size: $n\in\bN$

\State Create $(n\times(d+1))$ matrix $\mathbf{A}=(a_{ij})$

\State $i \gets 1$

\While{$i\leq n$}
\State Independently sample $Z_1,\dots,Z_{d+1}\sim N(0,1)$

\State $j\gets 1$

\While{$j\leq d+1$}

\State $a_{i j}\gets Z_j/\sqrt{Z_1^2+\dots+Z_{d+1}^2}$

\State $j\gets j+1$
\EndWhile
\State $i\gets i+1$
\EndWhile

\State \Return matrix $\mathbf{A}$, each row vector of which will follow a uniform distribution on $\bS^d$

\end{algorithmic}
\end{algorithm}

\begin{figure}
    \centering
    \includegraphics[width=1\linewidth]{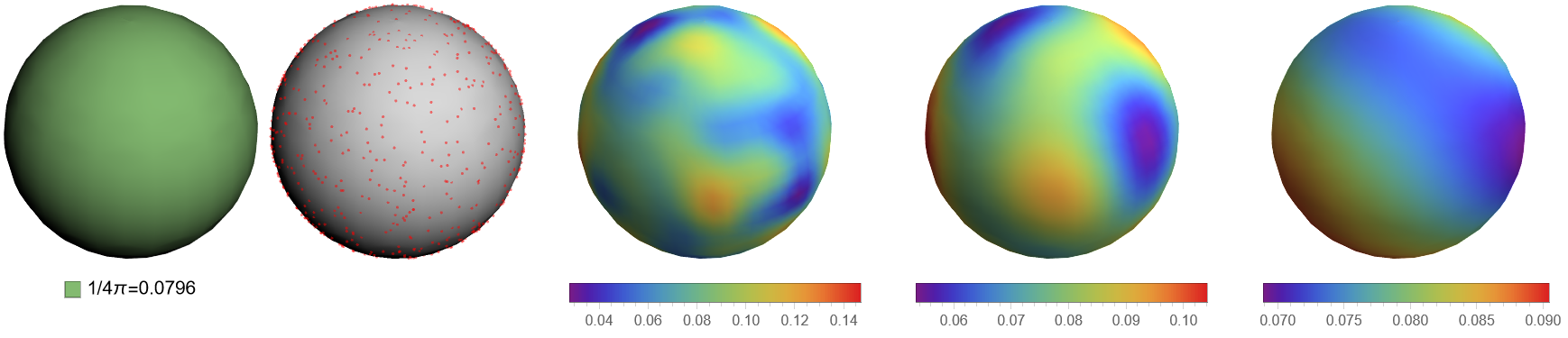}
    \caption{(Left to right) The plot of the uniform distribution over $\bS^2$; a sample of $n=1000$ points from such a distribution; KDEs constructed using said data for $s=0.5$, $s=1$ and $s=2$, respectively.}
    \label{fig:uniform sphere}
\end{figure}

We plot and evaluate our new KDEs in Figure \ref{fig:uniform sphere} and Table \ref{tab:unif probs}. 

The plots indicate that the estimators have captured the simulated data.

The KDE estimated probabilities (taken over halves and quarters of $\bS^2$) shown in Table \ref{tab:unif probs} are all comparable with both the observed frequencies and true probabilities, confirming the performance of the estimators introduced in the paper. 

Note finally that a preference to larger smoothness levels, say $s\in[1,2]$, seems better for such a uniform data, because it can avoid overfitting the data and ends up to very low MISE averages.

\begin{table}[]
    \centering
    \begin{tabular}{|c|c|c|c|c|c|}
    \hline
   Region & $s=0.5$&$s=1$&$s=2$& Frequency&True Prob.\\
     & $N_s=19$&$N_s=8$&$N_s=4$& & \\
    \hline
        \hline
        Half 1 $[0,\pi]\times(-\pi,0]$&$0.5100$&$0.5080$&$0.5077$&$0.518$&$0.5$\\
        Half 2 $[0,\pi]\times[0,\pi]$&$0.4900$&$0.4920$&$0.4923$&$0.482$&$0.5$\\
        \hline
        Quarter 1 $[0,\pi/2]\times(-\pi,0]$&$0.2454$&$0.2439$&$0.2422$&$0.246$&$0.25$\\
        Quarter 2 $[0,\pi/2]\times[0,\pi]$&$0.2471$&$0.2495$&$0.2464$&$0.243$&$0.25$\\
        Quarter 3 $[\pi/2,\pi]\times(-\pi,0]$&$0.2646$&$0.2641$&$0.2655$&$0.272$&$0.25$\\
        Quarter 4 $[\pi/2,\pi]\times[0,\pi]$&$0.2429$&$0.2425$&$0.2459$&$0.239$&$0.25$\\
        \hline
        \hline
        MISE & $0.00643$ &$ 0.00204$& $0.00062$& - &- \\
        \hline
    \end{tabular}
    \caption{Estimated probabilities from a uniform distribution on $\bS^2$ for $n=1000$ points and $s=0.5$, $s=1$, $s=2$, as well as the corresponding estimated MISE.}
    \label{tab:unif probs}
\end{table}

\begin{figure}
    \centering
    \includegraphics[width=1\linewidth]{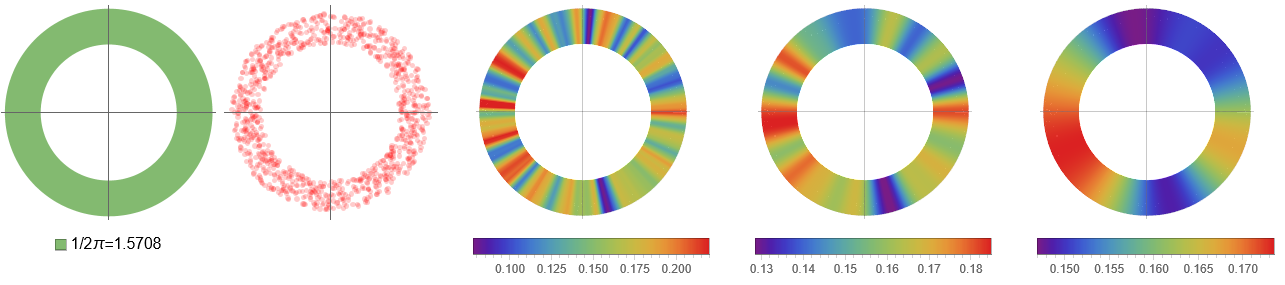}
    \caption{(Left to right) The plot of the uniform distribution over $\bS^1$; a sample of $n=1000$ points from such a distribution; KDEs constructed using said data for $s=0.5$, $s=1$, and $s=2$, respectively.}
    \label{fig:uniform circle}
\end{figure}

\paragraph{Scenario 2 (Uniform distribution on $\mathbb{S}^1$).}
Analogously, Figure \ref{fig:uniform circle} and Table \ref{tab:unif probs s1} contain the plots of the KDEs and the associated probabilities with the estimated MISE on $\bS^1$. Again the performance is satisfactory and $s\in[1,2]$ seems enough to avoid overfitting and offering low MISEs.

\begin{table}
    \centering
    \begin{tabular}{|c|c|c|c|c|c|}
    \hline
   Region & $s=0.5$&$s=1$&$s=2$& Frequency& True Prob.\\
     & $N_s=85$&$N_s=17$&$N_s=6$& &\\
    \hline
        \hline
        Half 1 $(-\pi,0]$ &$0.5124$&$0.5153$&$0.5130$&$0.515$&$0.5$\\
        Half 2 $[0,\pi]$&$0.4876$&$0.4847$&$0.4870$&$0.485$&$0.5$\\
        \hline
        Quarter 1 $(-\pi,-\pi/2]$&$0.2653$&$0.2676$&$0.2639$&$0.266$&$0.25$\\
        Quarter 2 $[-\pi/2,0]$&$0.2471$&$0.2477$&$0.2491$&$0.249$&$0.25$\\
        Quarter 3 $[0,\pi/2]$&$0.2372$&$0.2394$&$0.2387$&$0.233$&$0.25$\\
        Quarter 4 $[\pi/2,\pi]$&$0.2504$&$0.2453$&$0.2483$&$0.252$&$0.25$\\
        \hline
        \hline
        MISE & $0.00799$& $0.00248$&$0.00092$ &- &- \\
        \hline
    \end{tabular}
    \caption{Estimated probabilities from a uniform distribution on $\bS^1$ for $n=1000$ points and $s=0.5$, $s=1$, $s=2$, and corresponding estimated MISE.}
    \label{tab:unif probs s1}
\end{table}

\subsection{The von Mises-Fisher's distribution} 
We shall now move to estimating samples from the von Mises-Fisher (vMF) distribution; the spherical analogue of the normal distribution. 

Recall that the vMF distribution depends on two parameters: the mean direction $\mu\in\bS^d$ and the shape $\kappa>0$. Its density its given by
\begin{equation}
\label{vMF}
f(x)=c_{d,\kappa}\exp\big(\kappa \langle x,\mu\rangle\big)=c_{d,\kappa}e^{2\kappa}\exp\big(-\frac{\kappa}{2} \|x-\mu\|^2\big),\quad x\in\bS^d,
\end{equation}
where $c_{d,\kappa}$ is the normalisation constant. From \cite{MarJup}, we have that for $d=1$, this constant does not have a closed form and involves the modified Bessel function of the first kind: $c_{1,\kappa}=1/(2\pi I_0(\kappa))$. However, for $d=2$, one readily finds that the constant is  $c_{2,\kappa}=\kappa/(4\pi \sinh(\kappa))$.

The vMF density $f$ has a single peak at $\mu$ and the decay as $x$ departs from $\mu$ depends only on the distance between $x$ and $\mu$, see \eqref{metric on S2}, and the value of $\kappa$. Indeed, $\kappa$ plays the role of the reciprocal of the variance in a univariate Gaussian distribution with larger values of $\kappa$ corresponding to a greater concentration of the density about $\mu$. As $\kappa\to 0$, the vMF distribution approaches the spherical uniform distribution.

We now briefly present the methods we used to simulate from the vMF distribution. In the following, we assume the vMF distribution has a shape parameter $\kappa>0$ and mean direction $\mu=(1,0)$, in the $\bS^1$ case, or $\mu=(0,0,1)$, in the $\bS^2$ case. For other mean directions, the data need only be appropriately rotated after sampling. 

The algorithm, presented in \cite{BF1979}, for sampling over the vMF distribution on $\bS^1$, after fixing $\kappa>0$, is as given in Algorithm \ref{alg: vmf sim s1}.
\begin{algorithm}
    \caption{Simulation of von Mises-Fisher distributed points on $\bS^1$}\label{alg: vmf sim s1}
    \begin{algorithmic}[1]
        \Require Mean direction vector: $\mu=(\mu_1,\mu_2)\in\bS^1$
        \Require Concentration parameter: $\kappa>0$
        \Require Sample size: $n\in\bN$

        \State Create $(n\times 2)$ matrix $\mathbf{M}=(m_{ij})$
        \State $a \gets \sqrt{1+4 \kappa^2}$
        \State $b\gets (a-\sqrt{2a})/2\kappa$
        \State $r \gets (1+b^2)/2b$

        \State $i\gets 1$
        \While{$i \leq n$}
        \State $p\gets 0$
        \State $q \gets 0$
        \While{$p\leq 0$ and $q\leq 0$}
         \State Independently sample $U_1,U_2,U_3\sim$Unif$[0,1]$
        \State $z\gets \cos (\pi U_1)$
        \State $f \gets (1+rz)/(r+z)$
        \State $c\gets (r-f)\kappa$
        \State $p\gets c(2-c)U_2$
        \State $q \gets \ln(c/U_2)+1-c $
        \EndWhile 
       
        \State $\vartheta \gets$sign$(U_3-1/2)\arccos f$
\State $m_{i1}\gets \cos \vartheta$
\State $m_{i2}\gets \sin \vartheta$ \Comment{The rows of $\mathbf{M}$ will be vMF with mean $(1,0)$}
\State $i\gets i+1$
        \EndWhile
\If{$\mu_2\geq 0$}
\State $\phi \gets \arccos \mu_1$
\Else $\;\phi \gets -\arccos \mu_1$
\EndIf
\State $\mathbf{R}\gets\begin{pmatrix}
    \cos \phi & \sin \phi \\ -\sin \phi & \cos \phi
\end{pmatrix}$\Comment{Make a rotation matrix}
\State $\mathbf{V}\gets \mathbf{M} \,\mathbf{R}$

\Return matrix $\mathbf{V}$, the row vectors of which follow a von Mises Fisher distribution on $\bS^1$ with concentration parameter $\kappa$ and mean direction $\mu$.
    \end{algorithmic}
\end{algorithm}

In order to simulate from the von Mises-Fisher distribution on $\bS^2$, we follow the process given in \cite{Wood vmf}, which is given in Algorithm \ref{alg: vmf sim s2}.

\begin{algorithm}
    \caption{Simulation of von Mises-Fisher distributed points on $\bS^2$}\label{alg: vmf sim s2}
    \begin{algorithmic}[1]
    \Require Mean direction vector: $\mu=(\mu_1,\mu_2,\mu_3)\in \bS^2$
    \Require Concentration parameter: $\kappa>0$
    \Require Sample size: $n\in \bN$
    \State Create $(n\times 3)$ matrix $\mathbf{M}=(m_{ij})$
    \State $i\gets 1$
    \While{$i\leq n$}
    \State Sample $\varphi\sim$Unif$(-\pi,\pi]$
    \State Sample $U\sim$Unif$[0,1]$
    \State $\vartheta \gets \arccos (1+\frac{1}{\kappa}\ln(U+(1-U)e^{-2\kappa}))$
    \State $m_{i,1}\gets \cos \varphi \sin \vartheta$
    \State $m_{i,2}\gets \sin \varphi \sin \vartheta$
    \State $m_{i,3}\gets \cos \vartheta$ \Comment{The rows of $\mathbf{M}$ will be vMF with mean $(0,0,1)$}
\State $i \gets i+1$
\EndWhile
\If{$\mu\neq(0,0,-1)$}
\State $\mathbf{R}\gets \begin{pmatrix}
   1-\frac{\mu_1^2}{1+\mu_3}& -\frac{\mu_1\mu_2}{1+\mu_3}& -\mu_1 \\
   -\frac{\mu_1\mu_2}{1+\mu_3}& 1-\frac{\mu_2^2}{1+\mu_3}& -\mu_2 \\
   \mu_1 & \mu_2 & \mu_3
 \end{pmatrix}$ \Comment{Make a rotation matrix}
 \Else 
 
 $\mathbf{R}\gets - \mathbf{I}$
 \EndIf
 \State $\mathbf{V}\gets \mathbf{M}\, \mathbf{R}$

\Return matrix $\mathbf{V}$, the row vectors of which will be follow a von Mises-Fisher distribution with concentration parameter $\kappa$ and mean direction $\mu$.
    \end{algorithmic}
\end{algorithm}

\subsubsection*{Implementation}
\begin{figure}
    \centering
  \centerline{\includegraphics[width=1\linewidth]{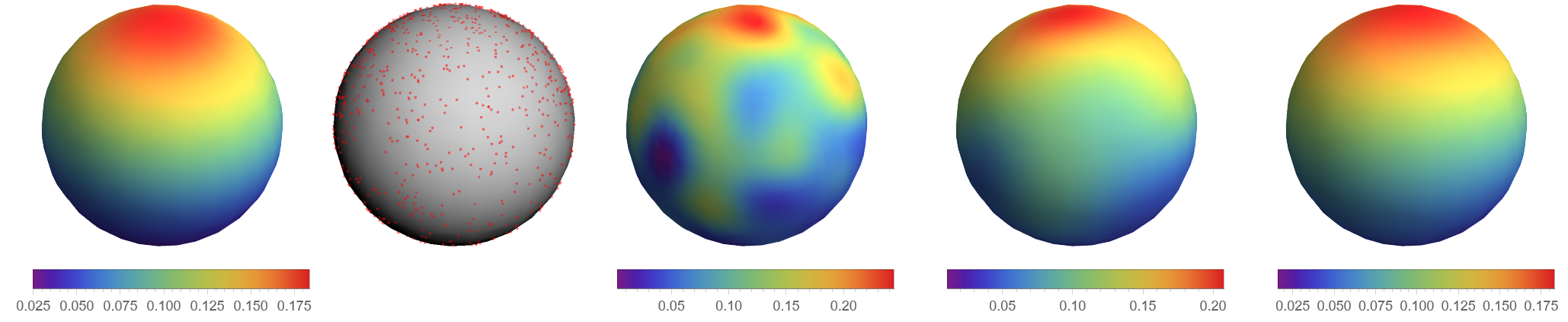}}
    \caption{(Left to right) The plot of the vMF distribution over $\bS^2$ with $\kappa=1$ and $\mu=(0,0,1)$; a sample of $n=1000$ points from such a distribution; KDEs constructed using said data for $s=0.5$, $s=1$ and $s=2$, respectively.}
    \label{fig:vmf k1}
\end{figure}

\paragraph*{Scenario 3 (vMF distribution on $\mathbb{S}^2$).}
In Figure \ref{fig:vmf k1} we plot the $\kappa=1$, $\mu=(0,0,1)$ vMF distribution alongside the $n=1000$ sample and three KDEs for $s=0.5$, $1$, $2$. Table \ref{tab:vmf k1 probs} contains the observed frequencies and true probabilities over four regions centred on the mean direction, alongside the MISE corresponding to each value of $s$.

\begin{table}[]
    \centering
    \begin{tabular}{|c|c|c|c|c|c|}
    \hline
    Region & $s=0.5$&$s=1$&$s=2$& Frequency&True Prob.\\
     & $N_s=19$&$N_s=8$&$N_s=4$& &\\
    \hline
    \hline
       $[0,\pi/2]\times(-\pi,\pi]$&0.7398&0.7414&0.7435&0.736&0.7311\\
       $[0,\pi/3]\times(-\pi,\pi]$&0.4578&0.4573&0.4576&0.45&0.4551\\
       $[0,\pi/4]\times(-\pi,\pi]$&0.2913&0.2964&0.2925&0.289&0.2936\\
       $[0,\pi/5]\times(-\pi,\pi]$&0.2028&0.2058&0.1991&0.199&0.2011\\
       \hline 
       \hline
       MISE & 0.00672 & $0.00203$ & $0.00063$&- & - \\
       \hline
    \end{tabular}
    \caption{Estimated probabilities from a vMF ($\kappa=1$, $\mu=(0,0,1)$) distribution for $n=1000$ points and $s=0.5$, $s=1$, $s=2$.}
    \label{tab:vmf k1 probs}
\end{table}

\begin{figure}
    \centering
    \includegraphics[width=0.5\linewidth]{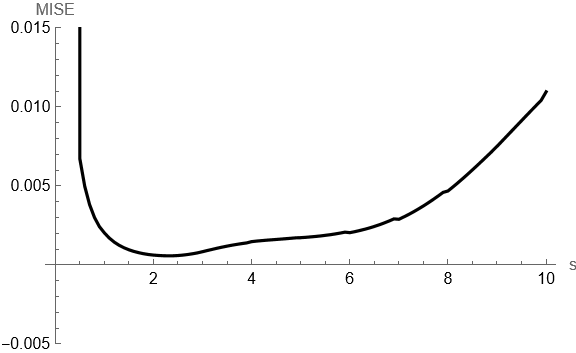}
    \caption{The estimated MISE for vMF simulated points with $\kappa=1$ and  $s\in\{0.1,0.2,\dots,10\}$, showing a clear minimum around $s=2$.}
    \label{fig: MISE for changing s}
\end{figure}

We can see that all values of $s$ gave very accurate estimations of the probabilities over each region. The plots and MISE both indicate that $s=2$ gives the best fit out of the values we considered. 
See Figure \ref{fig: MISE for changing s} for the MISE as $s$ ranges over $\{0.1,0.2,\dots,10\}$ where we see an approximate local minimum at $s=2.3$.  


\begin{figure}
    \centering
    \includegraphics[width=1\linewidth]{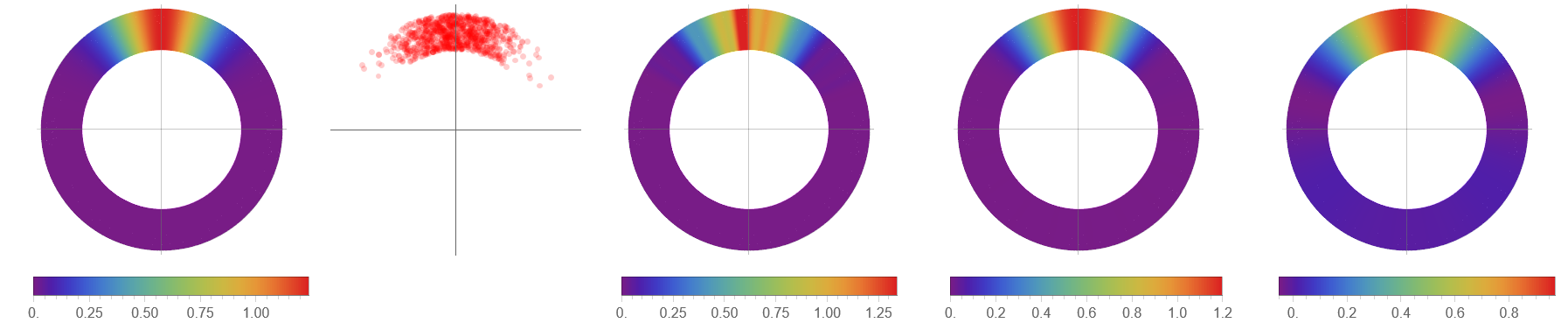}
    \caption{(Left to right) The plot of a vMF density over $\mathbb{S}^1$ with $\kappa=10$, $\mu=(0,1)$; a sample of $n=1000$ points from such a distribution; KDEs constructed using said data for $s=0.5$, $s=1$ and $s=2$, respectively.}
    \label{fig:vM S1}
\end{figure}

\paragraph*{Scenario 4 (vMF on $\mathbb{S}^1$).} 
We simulate $n=1000$ data from the vMF distribution on the circle with $\kappa=10$ and $\mu=(0,1)$. The findings are summarized in Figure \ref{fig:vM S1} and Table \ref{tab:vm k10 probs}.

\begin{table}[]
    \centering
    \begin{tabular}{|c|c|c|c|c|c|}
    \hline
    Region & $s=0.5$&$s=1$&$s=2$& Frequency&True Prob.\\
     & $N_s=85$&$N_s=17$&$N_s=6$& &\\
    \hline
    \hline
    $[0,\pi]$&0.999997&0.999944&0.999887& 1 & 0.999989\\
    $[\pi/3,2\pi/3]$&0.873325&0.874004& 0.817978&0.874&0.893612\\
    $[2\pi/5,3\pi/5]$&0.656978&0.649173&0.554804&0.662&0.671049\\
    $[5\pi/11,6\pi/11]$& 0.327823&0.328058&0.267136&0.329&0.343864\\

    \hline \hline 
    MISE &0.00762546&0.00214386&0.0445443  & -  & - \\ \hline
    \end{tabular}
    \caption{Estimated probabilities from a vMF density with $\kappa=10$ and $\mu=(0,1)$ on $\mathbb{S}^1$ for $n=1000$ points and $s=0.5$, $s=1$, $s=2$.}
    \label{tab:vm k10 probs}
\end{table}

\subsection{Mixture distributions}

Let us now consider von Mises-Fisher mixture densities. More precisely, consider some $m\in\bN$ and positive weights $w_1,\dots,w_m$ such that $\sum_j w_j=1$. Then for mean directions $\mu_1,\dots,\mu_m\in\bS^d$ and shape parameters $\kappa_1,\dots,\kappa_m>0$, we can define the vMF mixture density with $m$ components by
\begin{equation}
\label{eq: vMF mixture}
    f(x)=\sum_{j=1}^m w_j \cdot c_{d,\kappa_j} \exp{(\kappa_j \langle x,\mu_j\rangle)},\quad x\in\bS^d,\quad d=1,2,
\end{equation}
where $c_{d,\kappa_j}$ is the normalisation constant as before. This gives us a flexible family of distributions on the $d$-sphere $\bS^d$.  Let us denote by $W=(w_1,\dots,w_m)$, $K=(\kappa_1,\dots,\kappa_m)$, and $M=(\mu_1,\dots,\mu_m)$, the vectors of weights, shape parameters, and mean directions, respectively. 

\begin{figure}
    \centering
    \includegraphics[width=1\linewidth]{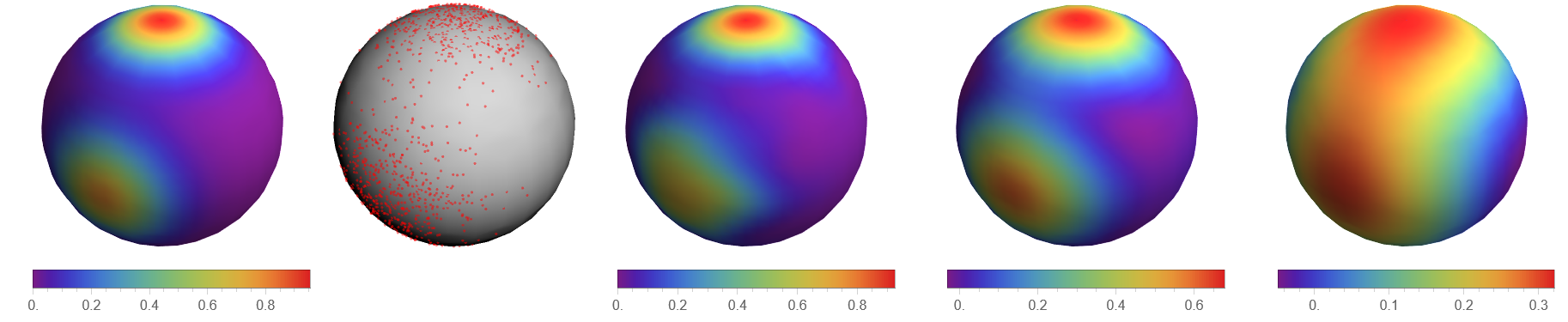}
    \caption{(Left to right) The plot of a mixture density over $\bS^2$; a sample of $n=1000$ points from such a distribution; KDEs constructed using said data for $s=0.5$, $s=1$ and $s=2$, respectively.}
    \label{fig:vMF mixed}
\end{figure}

\paragraph*{Scenario 5 (vMF mixture on $\mathbb{S}^2$).}
We plot such a density along with a sample of $n=1000$ points from it and KDEs constructed with $s=0.5$, $s=1$ and $s=2$ in Figure \ref{fig:vMF mixed}. The given density has $m=2$,  $W=(\frac{1}{2},\frac{1}{2})$, $K=(12,10)$ and $M=((0,0,1),(0,-1,0))$. Table \ref{tab:vmf mixed probs} contains the estimated probabilities for regions surrounding each of the two peaks, and corresponding estimated MISE.

\begin{table}[]
    \centering
    \begin{tabular}{|c|c|c|c|c|c|}
    \hline
    Region & $s=0.5$&$s=1$&$s=2$& Freq. &True Prob.\\
     & $N_s=19$&$N_s=8$&$N_s=4$& &\\
    \hline
    \hline
    $[0,\pi/6]\times(-\pi,\pi]$&0.3931&0.368&0.2324&0.394&0.4001\\
    $[0,\pi/4]\times(-\pi,\pi]$&0.4837&0.5024&0.4271&0.482&0.4886\\
        $[\pi/3,2\pi/3]\times[-2\pi/3,-\pi/3]$&0.4022&0.3854&0.2654&0.396&0.4012\\
        $[\pi/4,3\pi/4]\times[-3\pi/4,-\pi/4]$&0.4925&0.5021&0.466&0.492&0.4869\\
        \hline \hline
        MISE &0.0058 &0.017593& 0.15324 & - & -\\
        \hline
    \end{tabular}
    \caption{Estimated probabilities from a vMF mixture density on $\bS^2$ for $n=1000$ points and $s=0.5$, $s=1$, $s=2$.}
    \label{tab:vmf mixed probs}
\end{table}

Unlike Scenarios 1 and 3, both the MISE and plot show that $s=2$ is too large to provide a good fit in this case. This is expected as the density has more local behaviour since each mixture component vMF density has a very high shape parameter. Indeed, the $s=1$ plot and MISE indicate a better performance than $s=2$, with $s=0.5$ appearing as the best of the three. Moreover, the $s=2$ probabilities are much less accurate than the $s=0.5$ estimations, especially over smaller areas.

\paragraph*{Scenario 6 (vMF mixture on $\mathbb{S}^1$).}

\begin{figure}
    \centering
    \includegraphics[width=1\linewidth]{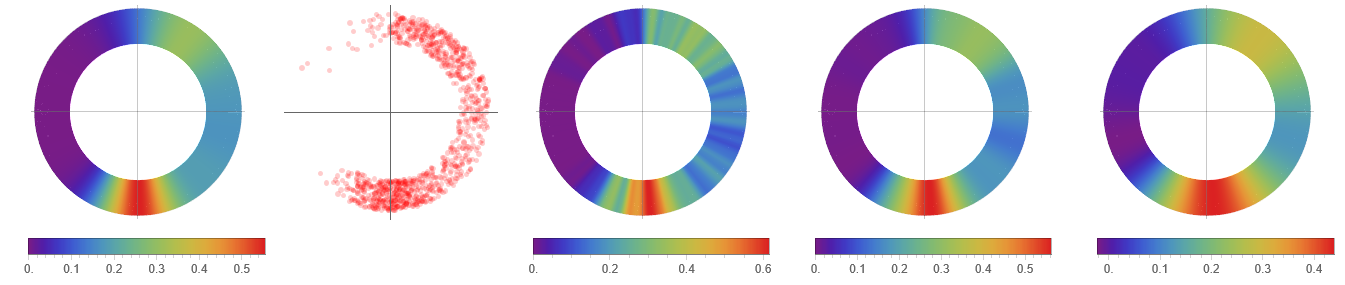}
    \caption{(Left to right) The plot of a mixture density over $\bS^1$; a sample of $n=1000$ points from such a distribution; KDEs constructed using said data for $s=0.5$, $s=1$ and $s=2$, respectively.}
    \label{fig:mixed S1}
\end{figure}

We conclude this section with the estimation of a vMF mixture density over on $\bS^1$. In this density we take $m=4$ components with weights $W=(\frac{1}{5},\frac{3}{10},\frac{1}{10},\frac{2}{5})$, shape parameters $K=(4,6,10,12)$, and mean directions $M=((1,0),(\frac{1}{2},\frac{\sqrt{3}}{2}),(\frac{1}{\sqrt{2}},\frac{1}{\sqrt{2}}),(0,-1))$. The corresponding plots and sample are given in Figure \ref{fig:mixed S1}, with probabilities and MISE given in Table \ref{tab:vm mixed probs}. 
Whilst the estimated probabilities are comparable for each value of $s$, the MISE and plots indicate that $s=1$ provides the best fit of the three. The $s=2$ plot shows oversmoothing and the $s=0.5$ plot shows undersmoothing.  

    \begin{table}[]
        \centering
        \begin{tabular}{|c|c|c|c|c|c|}
        \hline
        Region & $s=0.5$&$s=1$&$s=2$& Frequency&True Prob.\\
         & $N_s=85$&$N_s=17$&$N_s=6$& &\\
        \hline
        \hline
        $[-\pi/4,\pi/4]$&0.717&0.7213&0.7227&0.717&0.6968\\
        $[-3\pi/4,-\pi/4]$&0.5316&0.5285&0.5178&0.532&0.5394\\
        $[-\pi/2,0]$&0.6057&0.6091&0.615&0.602&0.5999\\
        $[\pi/12,7\pi/12]$&0.6469&0.647&0.6448&0.646&0.659\\
        \hline \hline 
        MISE & 0.00824& 0.00253& 0.01076 & -  & - \\ \hline
        \end{tabular}
        \caption{Estimated probabilities from a vMF mixture density on $\bS^1$ for $n=1000$ points and $s=0.5$, $s=1$, $s=2$.}
        \label{tab:vm mixed probs}
    \end{table}

\begin{remark}\label{Rem: small s}
We summarize our conclusions as follows:

(a) Near-uniform settings, such as Scenarios 1 and 2, as well as Scenario 3—where the vMF peak for $\kappa=1$ is very low—are well estimated using larger values of $s$ (like $s\approx2$), and equivalently a larger bandwidth $h$. 

(b) Moderate values of $s$ (such as $s\approx1$), and hence of $h$, are preferable when the underlying distribution departs from uniformity and exhibits more pronounced peaks along with flatter regions. 

(c) The only case where we observed a very small value of $s=0.5$  was Scenario 5, which involves two high peaks ($12$ and $10$) separated by a distance of $\pi/2$, leading to localized concentration in distinct regions with extended flat areas between them.
\end{remark}

\subsection{Comparison with vMF estimators}\label{sec: comparison 2}
Density estimation on the sphere dates back four decades. In particular, \cite{HWC} has attracted significant attention and represents one of the most widely used methods for handling spherical data. The kernel density estimators proposed in \cite[Section 5]{HWC} have the form
$$
\widetilde{f}_{\kappa}(x):=\frac{c_{d,\kappa}}{n}\sum_{i=1}^{n}\exp\big(\kappa\langle x,X_i\rangle\big),\quad \kappa>0,
$$
where $\kappa>0$ is a smoothing parameter playing the role of the reciprocal of a bandwidth, and $c_{d,\kappa}$ is as in (\ref{vMF}). In simple terms, the estimator $\widetilde{f}_{\kappa}$ averages vMF densities with common concentration parameter $\kappa$, centred at each observation $X_i$, and is often referred to as the vMF kernel density estimator (vMF-KDE).

In \cite[Section 5]{HWC}, the parameter $\kappa$ is selected via cross-validation methods by minimizing the squared-error loss and the Kullback--Leibler loss, respectively.
In this section, we compare our six simulation scenarios with the vMF-KDEs obtained via cross-validation methods in \cite{HWC}. Each of the two loss functions yields an estimated value of the parameter $\kappa$, which we denote by $\hat{\kappa}_{2}$ and $\hat{\kappa}_{\mathrm{KL}}$, respectively, and gives rise to the corresponding vMF-KDEs, denoted by $\widetilde{f}_{\hat{\kappa}_{2}}$ and $\widetilde{f}_{\hat{\kappa}_{\mathrm{KL}}}$.

The estimated $\mathrm{MISE}$ for each estimator across all scenarios is reported in Table \ref{tab:MISE of vMFKDEs}. As shown, the vMF-KDEs introduced in \cite{HWC} achieve smaller $\mathrm{MISE}$ values than the finite-order estimators of Theorem \ref{th: finite sum} in Scenarios 1 and 2, corresponding to the uniform distribution, when the target density is constant.

\begin{figure}
    \centering
    \includegraphics[width=0.9\linewidth]{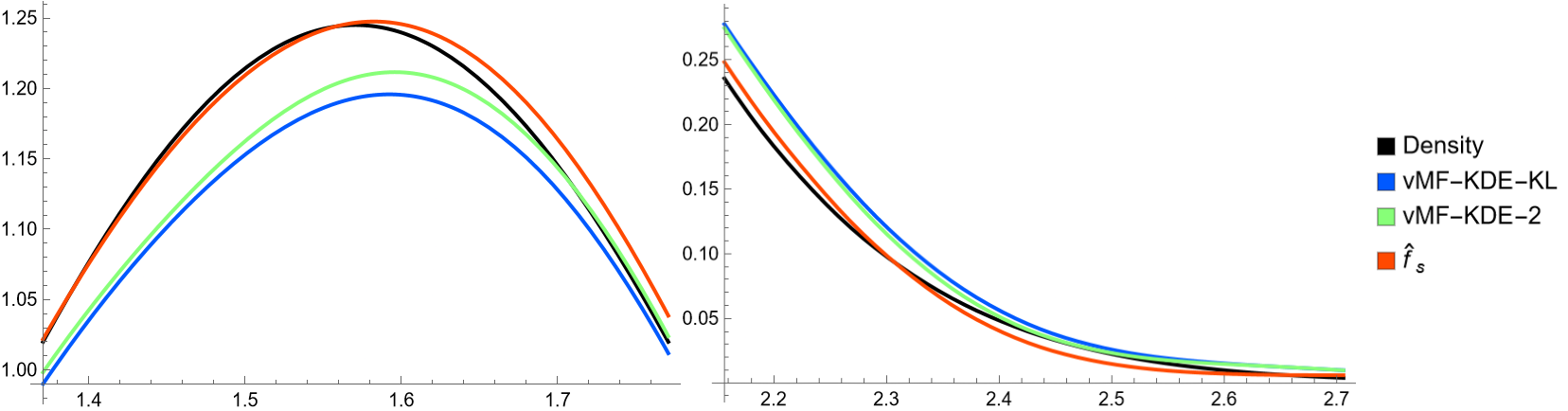}
    \caption{A closer view of the finite-order KDE in Scenario 5 in comparison to vMF-KDEs at the peak of the density (left) and the tail (right), plotted in terms of polar coordinates. The estimator $\widehat{f}_s$ captures both the peak and the tail.}
    \label{fig: KDE S1 close}
\end{figure}

In contrast, the finite-order estimators established here attain smaller $\mathrm{MISE}$ values in Scenarios 3--6, where the underlying densities (vMF distributions and their mixtures) exhibit one or more modes. This highlights the ability of our estimators to accurately capture local features of the density, including sharper extrema and appropriate decay in low-density regions; see also Figure \ref{fig: KDE S1 close}.

\begin{table}[]
\centering
\begin{tabular}{|c|c|c|c|}
\hline
Scenario & $\widehat{f}_s$ & $\widetilde{f}_{\hat{\kappa}_{2}}$ & $\widetilde{f}_{\hat{\kappa}_{\mathrm{KL}}}$ \\
\hline\hline
S1 & 
$(s=2);\,0.00062$ & 
$(\hat{\kappa}=0.7553)$; $0.000044$ & 
$(\hat{\kappa}=0.754);\,0.0000447$ \\
\hline
S2 & 
$(s=2);\,0.00092$ & 
$(\hat{\kappa}=2.0873);\,0.0001355$ & 
$(\hat{\kappa}=2.02467);\,0.000213$ \\
\hline
S3  & 
$(s=2);\,0.00063$ & 
$(\hat{\kappa}=11.237);\,0.0012$ & 
$(\hat{\kappa}=10.14);\,0.00115$ \\
\hline
S4 & 
$(s=1);\,0.00214$ & 
$(\hat{\kappa}=231.673);\,0.004822$ & 
$(\hat{\kappa}=150.204);\,0.0041813$ \\
\hline
S5 & 
$(s=0.5);\,0.0058$ & 
$(\hat{\kappa}=84.525);\,0.009502$ & 
$(\hat{\kappa}=63.5403);\,0.00961$ \\
\hline
S6 & 
$(s=1);\,0.00253$ & 
$(\hat{\kappa}=98.5197);\,0.00326289$ & 
$(\hat{\kappa}=62.09);\,0.00291$ \\
\hline
\end{tabular}
\caption{The estimated MISE of the KDE $\widehat{f}_s$, alongside the vMF-KDEs $\widetilde{f}_{\hat{\kappa}_{2}}$ and $\widetilde{f}_{\hat{\kappa}_{\mathrm{KL}}}$ for the value of $s>0$ that minimised MISE or the mean of the minimisers $\hat{\kappa}$ of each cross-validation function.}
    \label{tab:MISE of vMFKDEs}
\end{table}

\subsection{Comparison with adaptive density estimators}
Adaptive estimation is among the most important developments in nonparametric statistics and has been studied extensively over several decades across a variety of problems.

Adaptive density estimation on the sphere has been developed in \cite{BKMP,BKMP_Asymptotics} and further used in \cite{Geller2008}, through  wavelet estimators (note that wavelets are often referred to as needlets on compact manifolds), with coefficients selected via hard thresholding. This powerful methodology was introduced in \cite{DJKP} and later disseminated to a broader audience in the monograph \cite{HGPT}. In such adaptive estimators, there is no need to assume a given level $s$ of regularity, rather than a supremum value of the density $\|f\|_{\infty}$ and a parameter related with the threshold (parameter $k_0$ in \cite{BKMP}).

In \cite[Section 6]{BKMP}, adaptive wavelet estimators are numerically implemented for estimating a vMF mixture on $\mathbb{S}^2$, and compared with vMF-KDEs, using the (observed) $\mathcal{L}^{\infty}$-distance between the estimators and the target density. 

In this Section we compare our finite-order estimators with the adaptive wavelet estimators in the same scenario as in \cite{BKMP}: consider a mixture of $m=2$ vMF densities of the form \eqref{eq: vMF mixture}, with weights $W=(0.7,0.3)$, shape parameters $K=(1.2,8)$, and means $M=\big((0,1,0),\big(-\frac{\sqrt{2}}{2},-\frac{\sqrt{2}}{2},0\big)\big)$ and simulate $n=8000$ points from its density $f$.

\begin{figure}
    \centering
    \includegraphics[width=0.5\linewidth]{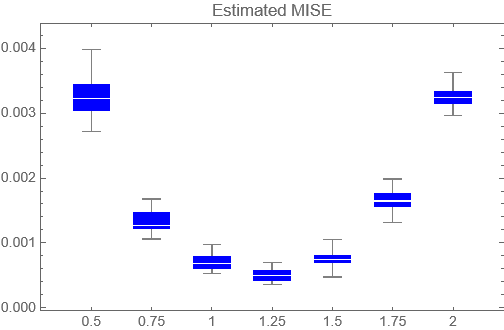}
    \caption{Boxplots for the estimated MISE of the finite-order KDE with $s\in\{0.5,0.75,1,1.25,1.5,1.75,2\}$, in the scenario as in \cite{BKMP}. The estimated MISEs are $0.00324$, $0.00133$, 0.00070, 0.00050, 0.00076, 0.00166, 0.00326, respectively.}
    \label{fig: MISE for changing s in the comparison with AoS}
\end{figure}

We estimate the $\mathcal{L}^\infty$-distance of our finite-order estimators, with $s=1$, following our Remark \ref{Rem: small s}, as the peaks are just moderately high. We refer further to the Figure \ref{fig: MISE for changing s in the comparison with AoS}, where one can see the estimated MISE for several values of $s$. The minimum estimated MISE came in the case of $s=1.25$ with comparable performance for $s=1$ and $s=1.5$.

More precisely, we evaluate the error of the estimator $|f(x)-\widehat{f}_s(x)|$ on a $0.01\times 0.01$ resolution grid over $[0,\pi]\times(-\pi,\pi]$, in spherical coordinates, giving approximately $2\times 10^5$ evaluations. For a given simulation, we record the maximum of the error $|f(x)-\widehat{f}_s(x)|$ over this grid. This is performed $30$ times and the $\mathcal{L}^\infty$-distance is estimated as the mean of these observed maximum errors. The results from this study, as well as the estimated $\mathcal{L}^\infty$-distance of the wavelet estimators for $k_0\in\{0.25,0.45,0.7\}$ are given in Table \ref{tab: adaptive Linf}. 

\begin{table}[]
    \centering
    \begin{tabular}{|c||c|}
    \hline
         Estimator& $\mathcal{L}^\infty$-distance  \\
         \hline \hline
        Wavelet, $k_0=0.25$ &  $0.073$\\
        Wavelet, $k_0=0.45$ & $0.045$\\
        Wavelet, $k_0=0.7$ & $0.066$\\
        $\widehat{f}_s$, $s=1$ & $0.0307$\\
        $\widehat{f}_s$, $s=1.25$ & $0.0293$\\
        \hline
    \end{tabular}
    \caption{The estimated $\mathcal{L}^\infty$-distance between the mixture density on $\mathbb{S}^2$ and the wavelet estimators and finite-order estimators, for $n=8000$ points.}
    \label{tab: adaptive Linf}
\end{table}

From this table we can see that our finite-order estimators succeed in achieving a lower $\mathcal{L}^\infty$-distance than the wavelet estimators. This is because our estimators are able to simultaneously capture peaks and tails; see the relevant discussion in \cite[Section 6]{BKMP}. 

In summary, our proposed estimators demonstrate that they are flexible enough to estimate multimodal densities with high peaks, a key advantage of adaptive wavelet density estimators, as noted in \cite{BKMP}. Of course, the value of any adaptive estimator is out of any question and is beyond any metrics. The main benefit of our finite-order estimators is their expression, which is more straight-forward to implement and does not require the structural elements behind the wavelet estimators of \cite{BKMP}, which are constructed via the theory in \cite{NPW}.

\subsection{Comparison of Theorem \ref{th: probability S2} with numerical integration}\label{sec: comparison}
Let us return to the motivational path to Theorem \ref{th: probability S2} and do the next simple experiment. 

\begin{figure}[h]
    \centering
    \includegraphics[width=1\linewidth]{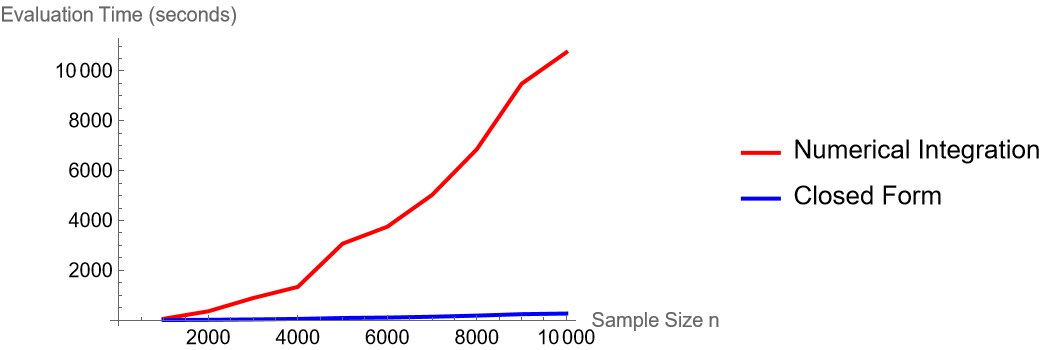}
    \caption{Evaluation times of the integral of the finite KDEs over a region $\bS^2$, using numerical integration (red) versus a closed form (blue) for $n\in\{1000,\dots,10\,000\}$.}
    \label{fig: s1 s2 integration times}
\end{figure}

Consider $n=1000,2000,\dots,10000$ data simulated uniformly on the sphere $\bS^2$.
The finite estimator \eqref{eq: kde S2 finite} was considered with a symbol $g_r(\lambda)$, as in Remark \ref{Remark 1}, where we took $s=1$ and $r=6$. The region of integration was $[\pi/4,3\pi/4]\times[-\pi/2,\pi/2]$.  We can see in Figure \ref{fig: s1 s2 integration times} that the closed form, established in Theorem \ref{th: probability S2}, significantly outperforms, in terms of time execution, the numerical integration method. For this simulation we used \textit{Mathematica}'s \texttt{GlobalAdaptive} method with the \texttt{GaussKronrodRule}, the default strategy of the native \texttt{NIntegrate} function.

\section{Case Studies}\label{sec: case studies} 
In this section we present four case studies, two on the unit circle $\bS^1$, and two on the unit sphere $\bS^2$. We illustrate how kernel density estimators on $\bS^1$ are well suited for spatially and temporally periodic data. In particular, we analyze the orientations of honeybees and the annual precipitation cycle in Los Angeles. For data on $\bS^2$, we construct KDEs for global datasets: earthquake locations and the directional distribution of stars in our galaxy. For each dataset, we compute a suitable finite–order KDE, as described in Theorem~\ref{th: finite sum}, and estimate probabilities over regions of interest using Proposition~\ref{th: integral S1 kde} for $\bS^1$ and Theorem~\ref{th: probability S2} for $\bS^2$.

\subsection{Spatially Periodic Data: Honeybee Orientations}
\label{sec: bees}
A natural source of periodic data is a random variable representing an angle or direction in physical space. Examples include wind and ocean current directions, rock core orientations, and animal movement headings. We illustrate this with the ``Honeybee Tracking I" dataset, made available by the Okinawa Institute of Science and Technology {\small\url{https://groups.oist.jp/bptu/honeybee-tracking-dataset}}. In this dataset, honeybees were filmed at 30 frames per second, and the orientation angle of each bee was recorded at every frame. Here, we analyze the first frame of ``Recording 1", originally gathered for the study \cite{BHMS}. Furthermore, the study \cite{BHPS} also gathered the trajectories of honeybees in a hive with a similar dataset available, also at the above URL.

\begin{figure}
\centering
\begin{subfigure}{.5\textwidth}
  \centering
  \includegraphics[width=.7\linewidth]{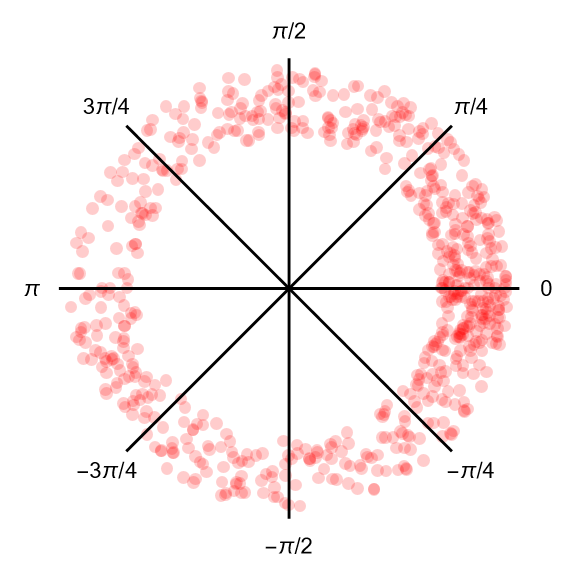}
\end{subfigure}%
\begin{subfigure}{.5\textwidth}
  \centering
  \includegraphics[width=.8\linewidth]{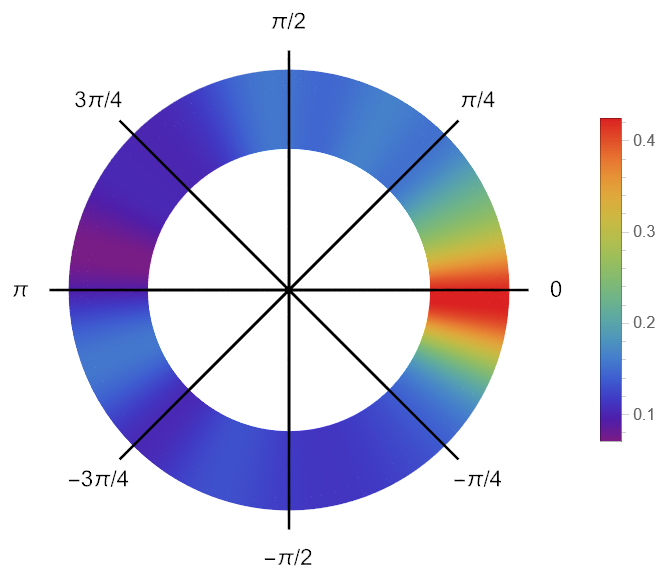}
\end{subfigure}
\caption{Plots of $n=691$ honeybee orientations (left) and the corresponding KDE with $s=1$ (right).}
\label{fig:Bees}
\end{figure}

In this dataset, bees are labelled as $0$, if their orientation is clear, or $1$, if their body is partially obscured within the hive. For consistency in our study, we only look at the $n=691$ honeybee orientations of class $0$, which are given in degrees. Converting these to radians makes these observations suitable for a KDE as in Theorem \ref{th: finite sum} and Remark \ref{Remark 1}. We selected $s=1$, and therefore $h=n^{-1/(2s+1)}\approx0.113$, $N_s=14$ and $r=5$. The data and corresponding KDE are plotted in Figure \ref{fig:Bees}. In this frame, the bees exhibit a clear preference for facing approximately towards the angle $-0.09$ radians with the least favoured position appearing to be approximately at 2.97 radians.

In Table \ref{tab: Bee probs} we give the estimated probabilities over two regions. The first is $[-1.33,0.81]$, which we estimate, from the KDE, to be the main `peak' of the density. The second region is $[2.46,\pi]\cup(-\pi,-3.03]$, which we estimate to be the `valley', or the region surrounding the global minimum. Note that for $[2.46,\pi]\cup(-\pi,-3.03]$, the estimated probability over this region is given by $\widehat\bP(A[2.36,\pi])+\widehat{\bP}(A(-\pi,-3.03])$, where $A[\vartheta_1,\vartheta_2]$ is given in \eqref{eq: region on S1}. We can see that the empirical and KDE estimated probabilities are approximately equal around the minimum and are relatively close around the maximum.

\begin{table}[]
    \centering
    \begin{tabular}{|c|c|c|}
    \hline
         Region& Frequency & KDE Prob. \\
         \hline
         \hline
       $[-1.33,0.81]$  & $0.4906$& $0.4893$\\
       \hline
       $[2.46,\pi]\cup(-\pi,-3.03]$ & $0.0637$&  $0.0637$   \\
       \hline
    \end{tabular}
    \caption{Empirical and KDE estimated probabilities of honeybee orientations around the maximum and minimum.}
    \label{tab: Bee probs}
\end{table}

Kernel density estimators on $\bS^1$ were particularly necessary in this case, rather than KDEs over $\bR$ as the periodicity of the data needed to be captured, both to properly represent the nature of directional data but also because in this case, the maximum occurred almost at $0$. 

\subsection{Temporally Periodic Data: Annual Precipitation Patterns}
We now turn to periodic precipitation patterns at a fixed location, illustrating our \textit{methodology} for treating temporal phenomena as circular data. The Open-Meteo API provides historical weather records for many locations worldwide; {\small\url{https://open-meteo.com/en/docs/historical-weather-api}}.

For Los Angeles, California, precipitation follows a clear annual cycle, with a wet season in spring and dry summers. From Open-Meteo, we extracted the days with non-zero precipitation from 01.01.2000 until 01.01.2025 (9133 observed days). This gave $n=1455$ days of precipitation in the time frame. Because the wet season occurs at approximately the same time each year, it is natural to view the data as periodic. 

More generally, if a random variable $T$ is periodic on $[a,b]$ with density $f_T:[a,b]\to[0,\infty)$ satisfying $f_T(a)=f_T(b)$, then independent observations $T_1,\dots,T_n\in(a,b]$ may be mapped to circular data on $\bS^1$ by
\begin{align*}
    [a,b] &\to [-\pi,\pi],\\
    t &\mapsto \frac{2\pi}{b-a}(t-a)-\pi.
\end{align*}
This transformation allows periodic temporal measurements to be analyzed using kernel density estimation on $\bS^1$, in line with our theory.

In our case, we assume that rainfall is periodic and therefore that the probability of rain occurring a second after 23:59:59 on the $31^\text{st}$ of December is the same as at 00:00:00 on the $1^\text{st}$ of January. We then can convert our observations of days of the year with rain occurring into angles on the circle $\bS^1$, via the above. Note that we adjust the period to be $366$ days in length on leap years.

\begin{figure}
\centering
\begin{subfigure}{.5\textwidth}
  \centering
  \includegraphics[width=.7\linewidth]{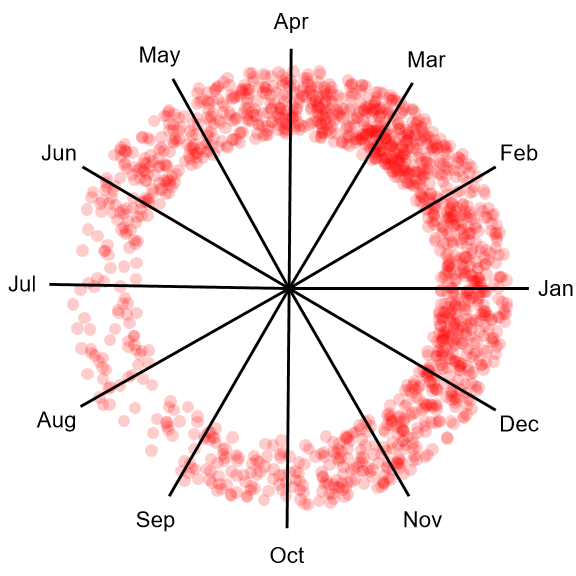}
\end{subfigure}%
\begin{subfigure}{.5\textwidth}
  \centering
  \includegraphics[width=.8\linewidth]{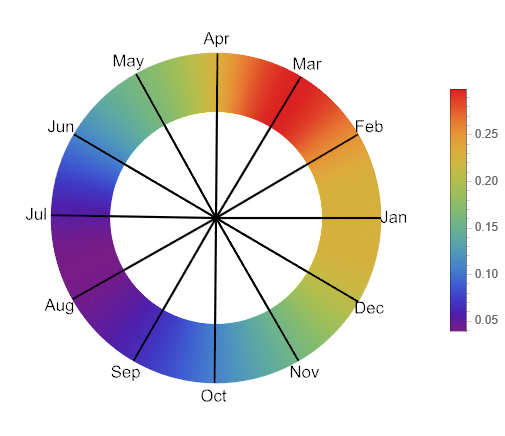}
\end{subfigure}
\caption{Plots of $n=1455$ days with nonzero precipitation in Los Angeles (left) and the corresponding KDE with $s=2$ (right).}
\label{fig: Rain LA}
\end{figure}

We applied the KDE from \eqref{eq: kde S1 finite} to the angular rainfall data $\theta_1,\dots,\theta_{1455}$, implementing the finite sum approximation from Theorem \ref{th: finite sum} with parameters $s=2$, $h\approx0.233$, $N_s=6$, and $r=6$. The resulting probability density, plotted in Figure \ref{fig: Rain LA} alongside the raw circular data, identifies the climatological regime of Los Angeles. A pronounced mode in the KDE during the spring accurately reflects the core of the wet season, while a broad minimum throughout the summer months captures the extended dry season. The first day of each month is marked on the plot, providing a calendar reference for this seasonal phenomenon.

In Table \ref{tab: LA rain probs} we give the estimated probabilities of rainfall over the period from the start of February to the end of March (the wet season), as well as the start of June to the end of October (the dry season). These KDE estimated probabilities are compared to the empirical probabilities in those regions. Both methods of probability estimation seem to agree in this table. We also performed this case study for $s=1$, finding that the estimated probabilities were approximately equal as to the $s=2$ case, with the difference being on the order of $10^{-3}$.

\begin{table}[]
    \centering
    \begin{tabular}{|c|c|c|}
    \hline
         Dates (Inclusive)& Frequency & KDE Prob.\\ 
         \hline
         \hline
       $1^\text{st}$ Feb. - $31^\text{st}$ Mar.  & $0.2838$& $0.2815$\\
       \hline
       $1^\text{st}$ Jun. - $31^\text{st}$ Oct. & $0.2027$ & $0.2008$\\
       \hline
    \end{tabular}
    \caption{Empirical and KDE estimated probabilities of rainfall in LA over certain periods.}
    \label{tab: LA rain probs}
\end{table}

By viewing seasonally dependent events as occurring on $\bS^1$, we are able to capture their periodic nature. If instead we were to view this data as having occurred over $\bR$ there would be two problems. Firstly, the kernel density estimator would estimate a very low probability for a date in  1990 or 2030, since it has no data for either period. Secondly, a particularly dry or wet year, i.e. an outlier, would greatly skew the estimation over that period of time, since there is at most one observation at any given day of the year. One could attempt to tackle both of these problems by using a KDE on a finite, non-periodic interval, say $[0,1]$ and scaling each observation accordingly. However, a KDE on this metric space wouldn't be periodic and there would likely be a discontinuity in the estimator as it moves from the 31st of December to the 1st of January. Therefore, the use of kernel density estimators over $\bS^1$ is essential in order to address the above problems and to properly construct an estimator of such a phenomenon.

\subsection{Global Data: Earthquake Locations}\label{sec: Seismology}

An obvious setting corresponding to the sphere $\bS^2$ is the surface of the Earth. With this in mind, we present an application of kernel density estimation to events on $\bS^2$ in estimating the distribution of earthquakes, as discussed in \cite{CGW}. Global earthquake data is freely available through the United States Geological Survey website \url{ https://earthquake.usgs.gov/earthquakes/search/}. 

There were $n=1630$ earthquakes of magnitude $6.5$ or higher from 1990 to 2024 (inclusive). Each earthquake recording comes with its latitude and longitude, which can easily be converted into Cartesian coordinates on $\bS^2$. Indeed, for a given latitude (lat) and longitude (lon) one can find the corresponding point on the sphere $X\in\bS^2$ via the following:
\begin{equation}
\label{eq: lat lon to cart}
    X=\Big(\cos\Big(\frac{\text{lat}\cdot\pi}{180}\Big)\cos\Big(\frac{\text{lon}\cdot\pi}{180}\Big),\cos\Big(\frac{\text{lat}\cdot\pi}{180}\Big)\sin\Big(\frac{\text{lon}\cdot\pi}{180}\Big),\sin\Big(\frac{\text{lat}\cdot\pi}{180}\Big)\Big).
\nonumber\end{equation}

\begin{figure}
\centering
\begin{subfigure}{.5\textwidth}
  \centering
  \includegraphics[width=1\linewidth]{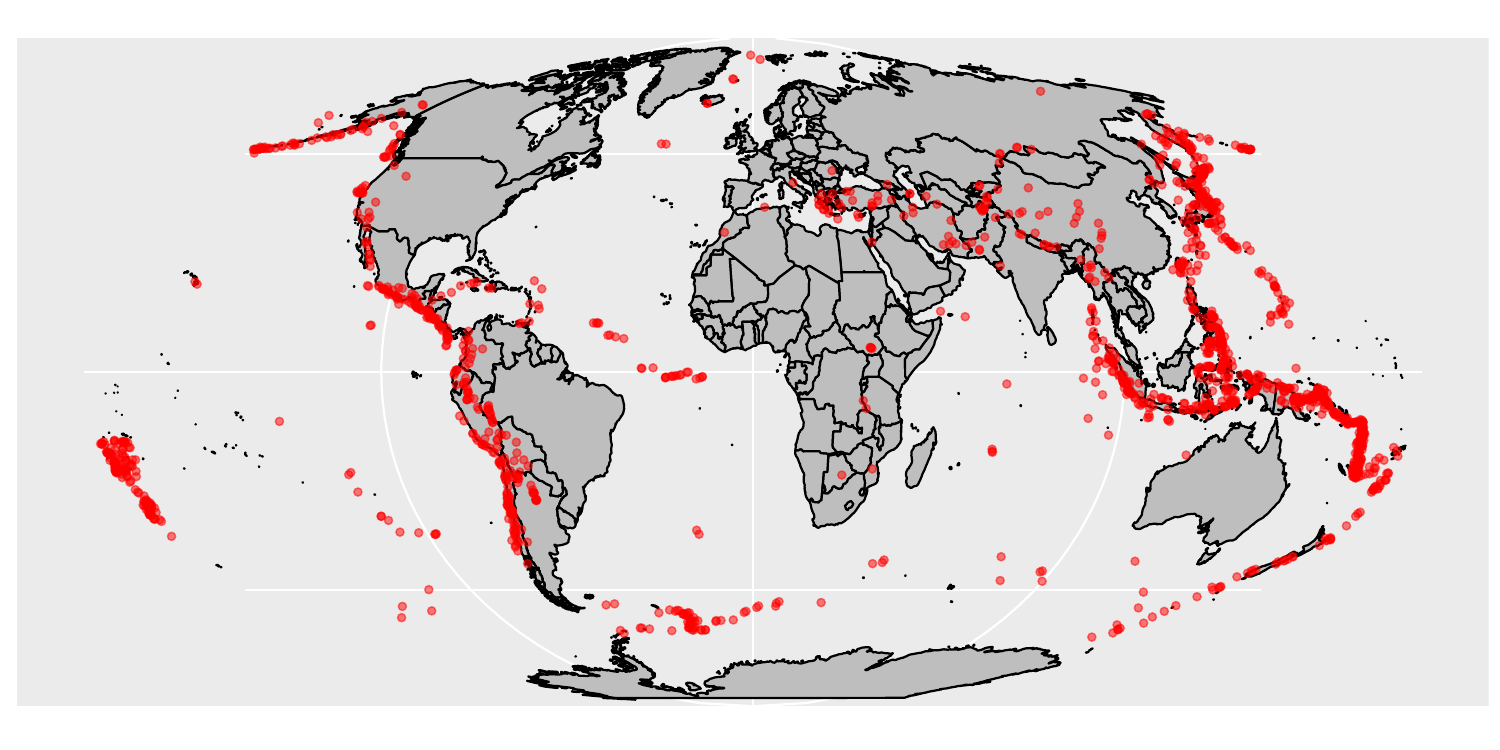}
\end{subfigure}%
\begin{subfigure}{.5\textwidth}
  \centering
  \includegraphics[width=1\linewidth]{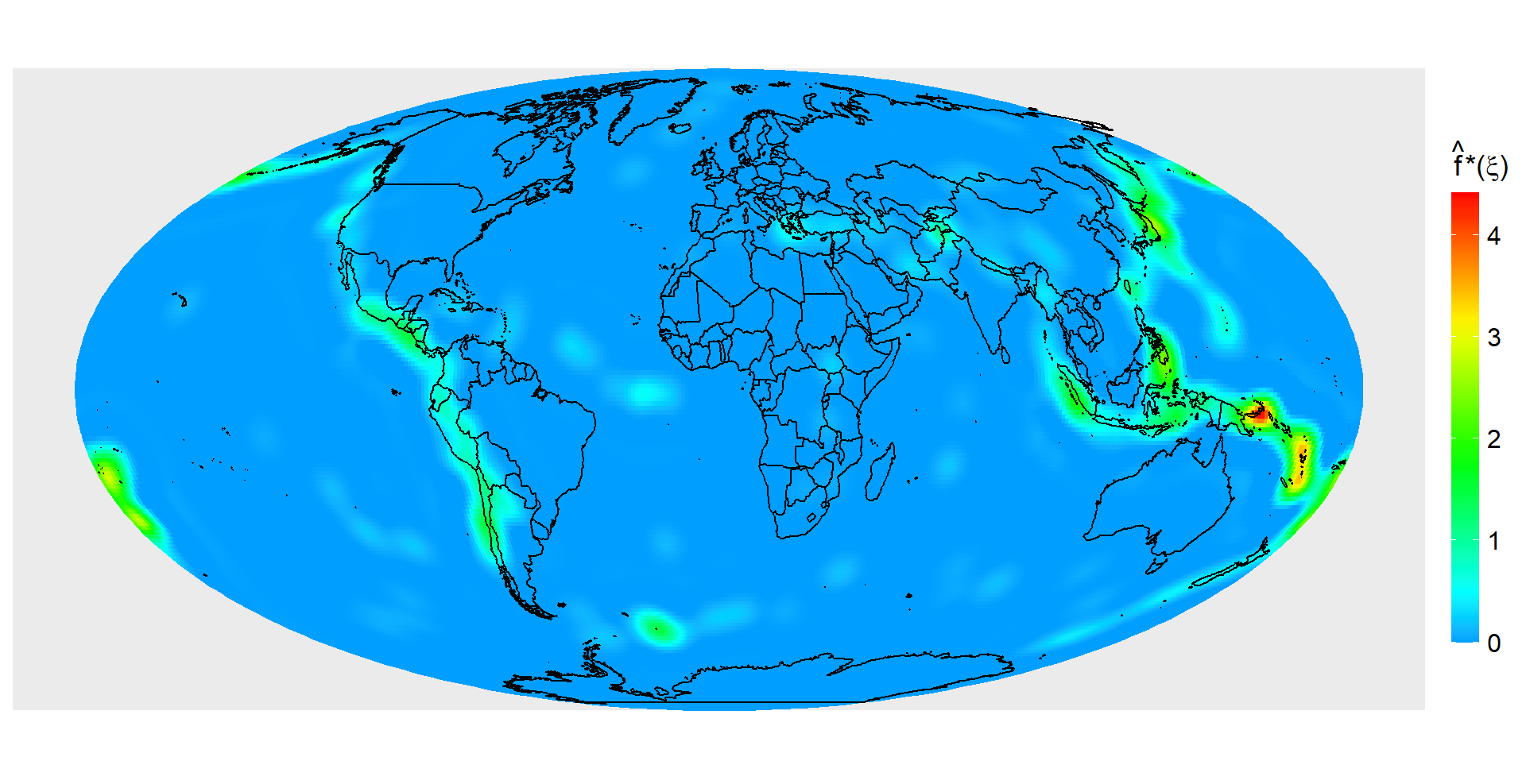}
\end{subfigure}
\caption{Plots of $n=1630$ high magnitude earthquakes (left) and the corresponding KDE with $s=0.05$ (right).}
\label{fig: earthquakes}
\end{figure}

With the above conversion in hand, we have observations $X_1,\dots,X_{1603}\in\bS^2$. We select $s=0.05$, as the phenomenon is clearly non-uniform and with several zones of high frequency; see Remark \ref{Rem: small s}. We construct a KDE $\widehat{f}_s$, as in Theorem \ref{th: finite sum} and Remark \ref{Remark 1}, taking $h=n^{-1/(2s+2)}\approx 0.0295$, $N_s=92$ and $r=6$. The corresponding KDE as well as the data is plotted in Figure \ref{fig: earthquakes}.

The resulting KDE effectively delineates regions of high seismic activity, with the Pacific 'Ring of Fire' appearing as the most prominent global feature. The model successfully captures the complex tectonics of Indonesia, where the Ring of Fire intersects with the Alpide Belt. We can use the probability estimation described in Section \ref{sec: probability} to estimate the probabilities over regions of particularly high (and low) seismic activity. We computed the probability of a high-magnitude earthquake within the bounding boxes of Japan, Chile, the Philippines, and Ireland by integrating the KDE over these regions. The results, presented in Table \ref{tab: earthquake probs} alongside the corresponding bounding coordinates and observed frequencies (to four decimal places), provide a probabilistic characterization of seismic hazard in these locations.

Note that when evaluating \eqref{eq: kde s2 integrated} for high cutoff points (about $N_s\geq60$, depending on the machine precision), it is recommended to using \textit{exact numbers}. The main reason for this is due to the numerical ill-conditioning of high-degree orthogonal polynomials, due to catastrophic cancellation. We found that using exact numbers was the simplest method to avoid this issue, however, there also exists numerically stable algorithms as alternatives. As an example, \textit{Mathematica} distinguishes between the \texttt{Real} (floating-point) number \texttt{0.9} and the \texttt{Rational} number \texttt{9/10}. The evaluation of an expansion of \texttt{LegendreP[60,x]} gives \texttt{-512} at \texttt{x=0.9}, clearly incorrect, whereas for \texttt{x=9/10} it gives \texttt{0.0317896}.

\begin{table}[]
    \centering
    \begin{tabular}{|c||c|c|c|c||c|c|}
    \hline
         \textit{Country} &\textit{Min. Lat.} &\textit{Max. Lat.} &\textit{Min. Lon.} &\textit{Max. Lon.}& Freq. &KDE Prob.  \\
         \hline
         \hline
      Japan   &$31^\circ$N&$45.5^\circ$N&$129.4^\circ$E &$145.5^\circ$E &$0.0540$&$0.049808$\\
      \hline
      Chile & $55.6^\circ$S&$17.6^\circ$S&$75.6^\circ$W&$70^\circ$W&$0.0429$&$ 0.043289$    \\
      \hline 
      Philippines & $5.6^\circ$N&$18.5^\circ$N&$117.2^\circ$E&$126.5^\circ$E&$0.0270$&$ 0.025305$\\
      \hline
      Ireland & $51.7^\circ$N&$55.1^\circ$N&$10^\circ $W&$6^\circ$W&$0.0000$&$ 0.000002$\\
      \hline
        \end{tabular}
    \caption{Empirical and KDE estimated probabilities of high magnitude earthquakes occurring within the bounding boxes of four regions on Earth.}
    \label{tab: earthquake probs}
\end{table}

This application highlight's the importance of having kernel density estimators available on $\bS^2$. If one were to have instead constructed an estimator over $\bR^2$, using the latitude and longitude as coordinates, the periodicity in the longitude would be lost. For instance, there would be a discontinuity between the earthquake density from Eastern Russia to Alaska. Furthermore, the distances between, say $(90^\circ \text{N},20^\circ \text{E})$ and $(90^\circ \text{N},30^\circ \text{E})$ would be taken by the KDE to be the same as the distance between $(0^{\circ}\text{N},20^{\circ}\text{E})$ and $(0^{\circ}\text{N},30^\circ\text{E})$, when in fact the former two points coincide. This would cause an underestimation of the density near the poles.

\subsection{Astronomical Data: Star Directions}\label{sec: Astronomy}
Distributions occurring on the sphere $\bS^2$ are not limited to just phenomena on Earth. Indeed, to each visible star, astronomers assign a ``galactic longitude" and ``galactic latitude". This galactic coordinate system assigns the sun of our solar system as its origin. The galactic equator is traced out from the distance of our sun to the approximate centre of the Milky Way. Using the right-hand convention, galactic longitude and galactic latitude are analogues of the familiar terrestrial coordinates and can therefore be converted to Cartesian coordinates via \eqref{eq: lat lon to cart}.

The Yale Bright Star Catalogue, freely available at \url{http://tdc-www.harvard.edu/catalogs/bsc5.html}, is a commonly used dataset of stars brighter than magnitude $6.5$. It contains the galactic longitude and latitude of $n=9096$ bright stars. As before, these can be easily converted into Cartesian coordinates $X_1,\dots,X_{9096}\in\bS^2$.

\begin{figure}
\centering
\begin{subfigure}{.5\textwidth}
  \centering
  \includegraphics[width=.9\linewidth]{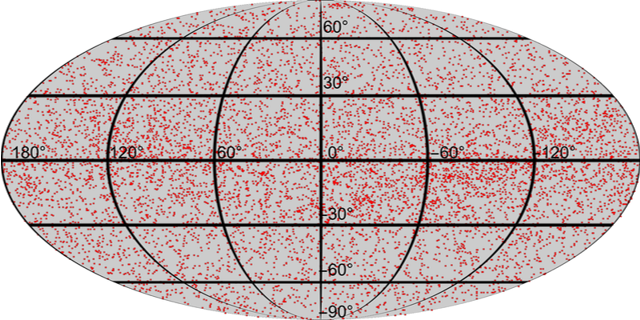}
\end{subfigure}%
\begin{subfigure}{.5\textwidth}
  \centering
  \includegraphics[width=.9\linewidth]{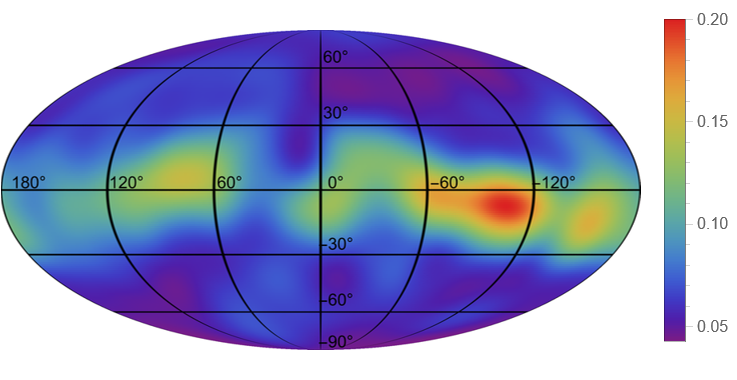}
\end{subfigure}
\caption{Plots of $n=9096$ bright stars (left) and the corresponding KDE with $s=1$ (right). }
\label{fig: stars}
\end{figure}

Figure \ref{fig: stars} presents the Bright Star Catalogue dataset alongside its (KDE) computed with $s=1$. This specific $s$ value was chosen to suppress small-scale noise while still resolving the major structural features of the distribution. The resulting KDE map, rendered in the standard galactic coordinate system where Earth is located at $(0^\circ,0^\circ)$ latitude and longitude, reveals regions of increased stellar density. In accordance with astronomical convention, the galactic longitude axis is reversed in the plots. A prominent band of high density is observed along the mid-plane of the Galaxy, an expected result that reflects the concentration of bright stars within the galactic disc. This concentration arises both from the higher density of stars in the disc and from the presence of luminous, massive stars which are largely confined to this plane.

In Table \ref{tab: star probs} we use Theorem \ref{th: probability S2} in order to estimate the probabilities over certain regions of space using the KDE. Indeed, we present the probabilities over the four quadrants of the galaxy, as given by galactic coordinates, as well as the probability over what appears to be the mode of the distribution, which we assume has a bounding box given by $(20^\circ$S, $130^\circ$W $)$ and $(5^\circ$N, $80^\circ$ W$)$. In the table we also present the empirical probabilities in each region. One can see that both the empirical probability and the KDE probability closely align. 

\begin{table}[]
    \centering
    \begin{tabular}{|c|c|c|c|c|c|c|c|}
    \hline
        Region & \textit{Min. Lat}&\textit{Max. Lat.}&\textit{Min. Lon.}&\textit{Max. Lon.}&Freq. &KDE Prob.\\
        \hline
        \hline
        Peak & $20^\circ$S&$5^\circ$N&$80^\circ$W&$130^\circ$W&$0.0629$&$0.0607$ \\
        \hline
        Upper Left & $0^\circ$N&$90^\circ$N&$180^\circ$W&$0^\circ$E&$0.233$&$0.2368$\\
        \hline
        Upper Right& $0^\circ$N&$90^\circ$N&$0^\circ$E&$180^\circ$E&$0.2379$&$0.2407$ \\
        \hline
        Lower Left &$90^\circ$S&$0^\circ$N&$180^\circ$W&$0^\circ$E&$0.2898$&$0.2847$\\
        \hline
        Lower Right &$90^\circ$S&$0^\circ$N&$0^\circ$E&$180^\circ$E&$0.2392$&$0.2379$\\
        \hline
    \end{tabular}
    \caption{The estimated (empirical and KDE) probabilities over certain regions of the sphere as well as the bounding coordinates of the regions, given in spherical coordinates.}
    \label{tab: star probs}
\end{table}

This case study demonstrates the usefulness of kernel density estimators over $\bS^2$ to estimate the distribution of the directions of vectors in $\bR^3$, where magnitude is irrelevant. 

\section{Conclusion}

In our study we provided nonparametric density and probability estimators on the circle $\bS^1$ and the sphere $\bS^2$ that are both rate-optimal and computationally efficient. The key idea is to replace the usual infinite expansions, as in \cite{CGKPP,CGP,CGW} underlying kernel density estimators with their finite-order counterparts, determining the order index in accordance with the smoothness of the density. The aforementioned works focused on the theoretical development of density estimation on broad metric spaces and on the further dissemination of the methods used therein. Our present study continues this line of work by focusing on computational aspects, which were not the main focus in \cite{CGKPP,CGP,CGW}. In particular, the need to implement these estimators in practical computations motivated us to investigate theoretically justified truncations of the infinite expansions. The finite-order estimators obtained in Theorem \ref{th: finite sum} preserve the optimality of the previous works and are ready to be used in data analysis. 

Moreover, our probability estimators admit closed-form expressions, which enables rapid computation even for large datasets. This strengthens the applicability of kernel methods in the analysis of circular and spherical data, where region-based probability summaries are often of primary interest.

We demonstrated a general methodology applicable to periodic and directional datasets and applied it in problems arising from zoology, climate science, seismology, and astronomy. These examples highlight the interpretability and computational advantages of working directly in the natural geometry of the data, rather than projecting it to Euclidean space. Our methodologies can be readily applied in a broad range of sciences.

There are several natural directions for future research. One is the extension of the finite-order methodology to higher-dimensional spheres $\bS^d$, $d>2$, or other metric spaces associated with Laplacians with discrete spectrum as in \cite{CGKPP}. Another natural direction is the use of cross-validation for data-driven selection of the bandwidth, which we plan to investigate in a future work.

\subsection*{Acknowledgements}
The authors would like to thank the Editor and the anonymous Referees for their careful reading and constructive comments, which have significantly improved the quality of this manuscript.

The first-named author would like to thank Prof.~Morgan Fraser for his helpful comments and suggestions on the astronomical application discussed in Section~\ref{sec: Astronomy}. 

The second-named author would like to thank Prof.~Dr.~Katarzyna Bozek an author of \cite{BHMS} and \cite{BHPS}, whose dataset and expertise enabled the discussion in zoology in Section~\ref{sec: bees}.

\section{Appendix}\label{sec: appendix} This section contains the proofs of our results.

\subsection{Proof of Theorem \ref{th: finite sum}} We start by proving the optimal rate of the finite estimators introduced in this study.

\begin{proof} The proof presents similarities, but also some differences for the case $d=1$ and $d=2$. Thus, we distinguish these cases.

\textit{Case } $d=2$.

Consider the MSE of our new estimators in (\ref{eq: kde S2 finite}) and use the trivial inequality $(a+b)^2\le 2(a^2+b^2)$;
\begin{align}\label{MSE proof 1}
\text{MSE}(\widehat f_s)(x)&
\leq 2\Big(\bE_f\big[\big(f(x)-\widehat f_{n,h}(x)\big)^2\big]+\bE_f\big[\big(\widehat f_{n,h}(x)-\widehat f_{s}(x)\big)^2\big]\Big)\\
&=2\Big(\text{MSE}(\widehat{f}_{n,h})(x) +\bE_f\big[\big(\widehat f_{n,h}(x)-\widehat f_{s}(x)\big)^2\big]\Big).
\nonumber\end{align}

From Theorem \ref{th: milk metrika} we know that MSE$(\widehat f_{n,h})(x)$ is bounded by $$c\, C(f)\,n^{-2/(2s+2)},$$ due to our selection of $h=n^{-1/(2s+2)}$. We now aim to bound the other term. By (\ref{eq: kde S2}) and (\ref{eq: kde S2 finite}) we arrive at

\begin{equation}
\big|\widehat{f}_{n,h}(x)-\widehat{f}_s(x)\big|=\Big|\frac{1}{n}\sum_{j=1}^{n}\sum_{\ell>N_s}\frac{ 2\ell+1}{4\pi}g(h\sqrt{\ell(\ell+1)})P_\ell(\langle x, X_j\rangle)\Big|.
\end{equation}
By (\ref{eq: decay}) we take
$$\big|g(h\sqrt{\ell(\ell+1)})\big|\le C_{\tau}\big(1+h\sqrt{\ell(\ell+1)}\big)^{-r}< C_{\tau}h^{-r}\big(\ell(\ell+1)\big)^{-r/2}.$$

For the Legendre polynomials we have $|P_{\ell}(u)|\le 1,$ for every $u\in[-1,1]$; see (\ref{legendre maximum}).

Then the error under consideration can be bounded from above by
\begin{align}\label{eq:trunc_error_sphere}
|\widehat{f}_{n,h}(x)-\widehat{f}_s(x)|&\le\frac{C_{\tau}}{4\pi}h^{-r}\sum_{\ell=N_s+1}^{\infty}(2\ell+1)\big(\ell(\ell+1)\big)^{-r/2}
\nonumber
\\
&\le \frac{C_{\tau}}{4\pi}h^{-r}\int_{N_s}^{\infty}(2x+1)\big(x(x+1)\big)^{-r/2}dx
\nonumber
\\
&=\frac{C_{\tau}}{2\pi(r-2)}h^{-r}\big(N_s(N_s+1)\big)^{1-\frac{r}{2}}
\nonumber
\\
&\le\frac{C_{\tau}}{2\pi(r-2)}h^{-r}N_s^{2-r}.
\end{align} 
Note that above we used the usual bound of the series of decreasing sequences by the corresponding integral and then we changed the variable of integration.

For the choice of the bandwidth it holds $h^{-r}=n^{r/(2s+2)}$. 

By the selected $N_s$ as in \eqref{eq: terminal N Sd} we get
$$N_s^{2-r}\le 2\pi(r-2)n^{-(s+r)/(2s+2)}.$$

Replacing in (\ref{eq:trunc_error_sphere}) we derive the bound
\begin{equation}\label{eq:boundofthefiniteerror}
|\widehat{f}_{n,h}(x)-\widehat{f}_s(x)|\le C_{\tau}n^{-s/(2s+2)}.
\end{equation}

Combining (\ref{MSE proof 1}) with (\ref{eq:boundofthefiniteerror}) and the result from Theorem \ref{th: milk metrika} we arrive at \eqref{eq: s2 trunc thm bound} and complete the proof for $d=2$.

\textit{Case } $d=1$. As in the previous case we focus in $|\widehat f_{n,h}(x)-\widehat f_s(x)|$.

For a given $\ell\in\bN$, we have using (\ref{eq: decay}) that 
\begin{equation}
|g(h\ell)|\le C_{\tau}(1+h\ell)^{-r}<(h\ell)^{-r}.
\end{equation}
Moreover, $|\cos(t)|\leq1$ for all $t\in\bR$. Therefore, we have the following bound:

\begin{align*}
|\widehat f_{n,h}(x)-\widehat f_s(x)|&\leq \frac{1}{\pi n}\sum_{j=1}^n\sum_{\ell=N_s+1}^\infty h^{-r}\ell^{-r}=\frac{h^{-r}}{\pi}\sum_{\ell=N_s+1}^\infty \ell^{-r}
\\
&\leq \frac{h^{-r}}{\pi}\int_{N_s} ^\infty u^{-r}du
=\frac{h^{-r}}{\pi}\frac{N_s^{1-r}}{r-1},
\nonumber\end{align*}
for every $x\in\bS^1$.

Now consider $N_s$ as in \eqref{eq: terminal N Sd}. One has that
\begin{equation*}
N_s^{1-r}\leq \pi(r-1)n^{-(s+r)/(2s+1)},
\end{equation*}
and substituting $h=n^{-1/(2s+1)}$, we get
\begin{equation*}
|\widehat f_{n,h}(x)-\widehat f_s(x)|\leq n^{-s/(2s+1)}.
\end{equation*}
and therefore, we have that
\begin{equation*}
\bE_f\Big[\big(\widehat f_{n,h}(x)-\widehat f_s(x)\big)^2\Big]\leq n^{-2s/(2s+1)}.
\end{equation*}
Combining this with (\ref{MSE proof 1}) and our bound for the MSE of $\widehat f_{n,h}$ from Theorem \ref{th: milk metrika} gives the result for $d=1$ and completes the proof.      
\end{proof}

\subsection{Proofs of the results in probability estimation} We present here all proofs of the results in Section \ref{sec: probability}. We start as follows:

\begin{proof}(of Proposition \ref{Th:MSEboundProbabilityEstimator}).
The probability of the measurable set $A$ is
\begin{equation*}
    \bP(A)=\int_{A}f(x)d\sigma_{\bS^d}(x).
\end{equation*}

Then by Jensen's inequality, for the uniform probability on $A$, and Fubini-Tonelli Theorem we have that
\begin{align*}
    \mathbb{E}_f&\big(\widehat{\bP}(A)-\bP(A)\big)^2\le
    \mathbb{E}_f\Big(\sigma_{\bS^d}(A)\int_{A}\big|\widehat{f}_s(x)-f(x)\big|\frac{d\sigma_{\bS^d}(x)}{\sigma_{\bS^d}(A)}\Big)^2 
    \\
    &\le
  \sigma_{\bS^d}(A) \bE_f  \int_{A}\big(\widehat{f}_s(x)-f(x)\big)^2 d\sigma_{\bS^d}(x)
    \\
    &= \sigma_{\bS^d}(A)  \int_{A}\bE_f \big(\widehat{f}_s(x)-f(x)\big)^2 d\sigma_{\bS^d}(x)\le \sigma_{\bS^d}(A)^2\sup_{x\in\bS^d}\text{MSE}(\widehat{f}_s(x)),
\end{align*}
and the result follows from Theorem \ref{th: finite sum}.

\end{proof}

The main probability estimation on $\bS^1$; Proposition \ref{th: integral S1 kde} has the following proof:
\begin{proof}
We are concerned with evaluating $$\int_{A[\vartheta_1,\vartheta_2]}\widehat f(x)d\sigma_{\bS^1}(x)=\int_{\vartheta_1}^{\vartheta_2} \widehat f_s(\theta) d\theta$$ where $\widehat f_s$ is as in \eqref{eq: kde S1 finite}. By simple integration we confirm the statement: 
\begin{align*}
\widehat{\bP}&(A[\vartheta_1,\vartheta_2])=\int_{\vartheta_1}^{\vartheta_2}\frac{1}{2\pi n}\sum_{j=1}^n\Big(1+2\sum_{\ell=1}^{N_s} g(h \ell)\cos(\ell(\theta-\theta_j))\Big)d\theta\nonumber\\
&=\frac{\vartheta_2-\vartheta_1}{2\pi }+\frac{1}{\pi n}\sum_{j=1}^n\sum_{\ell=1}^{N_s} g(h\ell)\int_{\vartheta_1}^{\vartheta_2}\cos(\ell(\theta-\theta_j))d\theta.
\nonumber\\    &=\frac{\vartheta_2-\vartheta_1}{2\pi }+\frac{1}{\pi n}\sum_{j=1}^n\sum_{\ell=1}^{N_s} g(h\ell)\frac{1}{\ell}(\sin(\ell(\vartheta_2-\theta_j))-\sin(\ell(\vartheta_1-\theta_j))).
\end{align*}
\end{proof}

Let us now proceed to prove the Lemma \ref{Lemma: integrals related to incomplete Beta function} on the integrals of the associated Legendre functions.
\begin{proof} 
Recall

First we combine (\ref{associated Legendre functions}) and (\ref{Simple for Legendre polynomials}). For every $\ell\in\bN$ and $1\le m\le\ell$, we differentiate (\ref{Simple for Legendre polynomials}) to extract
\begin{align}
\frac{d^m}{du^m}P_{\ell}(u)&=\sum_{k=m}^{\ell}{\scriptstyle\binom{\ell}{k}\binom{\ell+k}{k}}\frac{(-1)^k}{2^k}\frac{d^m}{du^m}(1-u)^k
\nonumber
\\
&=\sum_{k=m}^{\ell}{\scriptstyle \binom{\ell}{k}\binom{\ell+k}{k}}\frac{(-1)^{k+m}}{2^k}k(k-1)\cdots(k-m+1)(1-u)^{k-m}
\nonumber
\\
&=\sum_{k=m}^{\ell}{\scriptstyle \binom{\ell}{k}\binom{\ell+k}{k}}\frac{(-1)^{k+m}}{2^k}\frac{k!}{(k-m)!}(1-u)^{k-m}.
\end{align}
We embed it into (\ref{associated Legendre functions}) and obtain for every $\ell\ge 1,\;1\le m\le \ell$
\begin{align}
\label{final computation of associated Legendre functions}
P_{\ell}^{m}(u)&=
(-1)^m(1-u^2)^{m/2}\frac{d^m}{du^m}P_{\ell}(u)
\\
&=\sum_{k=m}^{\ell}{\scriptstyle \binom{\ell}{k}\binom{\ell+k}{k}}\frac{(-1)^{k}}{2^k}\frac{k!}{(k-m)!}(1+u)^{m/2}(1-u)^{k-m/2}.
\nonumber\end{align}

Therefore for the integral $\int_{\gamma_1}^{\gamma_2}P_{\ell}^m(u)du$, it suffices to compute $$\int_{\gamma_1}^{\gamma_2}(1+u)^{m/2}(1-u)^{k-m/2}du.$$ We consider the more general integral
\begin{equation}
O(a,b,\gamma_1,\gamma_2):=\int_{\gamma_1}^{\gamma_2}(1+u)^{a}(1-u)^{b}du,\quad a,b>-1.
\end{equation}
We change to the variable $t:=\frac{1+u}{2}$ and for $-1\le\gamma_1\le u\le\gamma_2\le 1$, we have that $0\le\frac{1+\gamma_1}{2}\le t\le \frac{1+\gamma_2}{2}\le 1$ and therefore
\begin{align}
O(a,&b,\gamma_1,\gamma_2)=2^{a+b+1}\int_{\frac{1+\gamma_1}{2}}^{\frac{1+\gamma_2}{2}}t^a(1-t)^bdu
\nonumber\\
&=2^{a+b+1}\Big[B\Big(\frac{1+\gamma_2}{2};a+1,b+1\Big)-B\Big(\frac{1+\gamma_1}{2};a+1,b+1\Big)\Big].
\label{O integral}
\end{align}
We can now apply (\ref{O integral}) to (\ref{final computation of associated Legendre functions}) for $a=m/2$ and $b=k-m/2$ to conclude the proof of the Lemma.
\end{proof}

\subsubsection*{Proof of Theorem \ref{th: probability S2}} We can now proceed presenting the proof of the main result in probability estimation.

\begin{proof}
We consider the new optimal estimator in (\ref{eq: kde S2 finite}) and we use the addition formula for Legendre polynomials and spherical harmonics, given in \eqref{Legendreand spherical harmonics}, taking

\begin{align}
\widehat f_s(x)&=\frac{1}{n}\sum_{j=1}^n\sum_{\ell=0}^{N_s} \frac{2\ell+1}{4\pi}g\Big(h\sqrt{\ell(\ell+1)}\Big)P_\ell(\langle x, X_j\rangle)
\nonumber\\
&=\frac{1}{n}\sum_{j=1}^n\sum_{\ell=0}^{N_s} g\Big(h\sqrt{\ell(\ell+1)}\Big)\sum_{m=-\ell}^{\ell}Y^m_{\ell}(\theta,\phi)\overline{Y_{\ell}^m(\theta_j,\phi_j)},
\label{proof prob 1}\end{align}
where $(\theta,\phi)$ the spherical coordinates of $x\in\bS^2$.

Based on the properties of complex numbers and the complex conjugation between spherical harmonics in (\ref{spherical harmonics conjugate}) we infer that
\begin{align}
\label{Spherical simplified}
\sum_{m=-\ell}^{\ell}&Y^m_{\ell}(\theta,\phi)\overline{Y_{\ell}^m(\theta_j,\phi_j)}\nonumber\\&=\sum_{m=-\ell}^{-1}Y^m_{\ell}(\theta,\phi)\overline{Y_{\ell}^m(\theta_j,\phi_j)}+Y^0_{\ell}(\theta,\phi)\overline{Y_{\ell}^0(\theta_j,\phi_j)}
\nonumber
\\&\hspace{1cm}+\sum^{\ell}_{m=1}Y^m_{\ell}(\theta,\phi)\overline{Y_{\ell}^m(\theta_j,\phi_j)}
\nonumber\\
&=Y^0_{\ell}(\theta,\phi)\overline{Y_{\ell}^0(\theta_j,\phi_j)}+2\sum^{\ell}_{m=1}\text{Re}\big(Y^m_{\ell}(\theta,\phi)\overline{Y_{\ell}^m(\theta_j,\phi_j)}\big).
\end{align}
Observe that for $m=0$, the corresponding spherical harmonics $Y_{\ell}^0$ are real-valued and precisely
\begin{equation}
Y^0_{\ell}(\theta,\phi)\overline{Y_{\ell}^0(\theta_j,\phi_j)}=N_{\ell,0}^2P_{\ell}(\cos\theta)P_{\ell}(\cos\theta_j).
\end{equation}
On the other hand, by the formula (\ref{spherical by associated Legendre}) we extract

\begin{align}
\nonumber\text{Re}&\big(Y^m_{\ell}(\theta,\phi)\overline{Y_{\ell}^m(\theta_j,\phi_j)}\big)=N_{\ell,m}^2P_{\ell}^m(\cos\theta)P^m_{\ell}(\cos\theta_j)\text{Re}\Big(e^{im(\phi-\phi_j)}\Big)
\\
&=N_{\ell,m}^2P_{\ell}^m(\cos\theta)P^m_{\ell}(\cos\theta_j)\cos\big(m(\phi-\phi_j)\big).
\label{proof prob 2}
\end{align}

Combining (\ref{proof prob 1})-(\ref{proof prob 2}), the estimated probability takes the form
\begin{align}
\widehat{\bP}(\cA_{\vartheta,\varphi})&=\int_{\cA_{\vartheta,\varphi}}\widehat{f}_{s}(x)d\sigma_{\bS^d}(x)\nonumber
\\
&=\frac{1}{n}\sum_{j=1}^n\sum_{\ell=0}^{N_s} g\Big(h\sqrt{\ell(\ell+1)}\Big)\times
\nonumber\\
&\hspace{1cm}\Bigg[\int_{\vartheta_1}^{\vartheta_2}\int_{\varphi_1}^{\varphi_2}N_{\ell,0}^2P_{\ell}(\cos\theta)P_{\ell}(\cos\theta_j)\sin\theta d\phi d\theta
\nonumber\\
&\hspace{1cm}+2\sum_{m=1}^{\ell}\int_{\vartheta_1}^{\vartheta_2}\int_{\varphi_1}^{\varphi_2}N_{\ell,m}^2P_{\ell}^m(\cos\theta)P^m_{\ell}(\cos\theta_j)\times
\nonumber
\\
&\hspace{4cm}\times\cos\big(m(\phi-\phi_j)\big)\sin\theta d\phi d\theta\Bigg]
\nonumber\\
&=:\frac{1}{n}\sum_{j=1}^n\sum_{\ell=0}^{N_s} g\Big(h\sqrt{\ell(\ell+1)}\Big)\Big(I_{\ell}^{0}(j)+2\sum_{m=1}^{\ell}I_{\ell}^{m}(j)\Big),
\end{align}
and we must prove that the integrals $I_{\ell}^{m}(j)$, agree with their values in the statement of the theorem.

In the integrals $I_\ell^{0}(j)$ first. We observe by changing the variable $u=\cos\theta$ that

\begin{equation}
I_{\ell}^{0}(j)=N_{\ell,0}^2(\varphi_2-\varphi_1)P_{\ell}(\cos\theta_j)\int^{\cos\vartheta_1}_{\cos\vartheta_2}P_{\ell}(u)d u.
\end{equation}

Recall the expansion (\ref{Simple for Legendre polynomials}) of the Legendre polynomials and focus on the integrals $\int^{\cos\vartheta_1}_{\cos\vartheta_2}(1-u)^k du$, $\le k\le\ell$.

It is easily observed (we may distinguish the cases $k=0$ and $k>0$ to compute the unified result), that for all $0\le k\le \ell$
\begin{equation*}
\int^{\cos\vartheta_1}_{\cos\vartheta_2}(1-u)^kdu=\frac{1}{k+1}\Big( (1-\cos\vartheta_2)^{k+1}-(1-\cos\vartheta_1)^{k+1}\Big),\quad k\ge0,
\end{equation*}
which concludes the calculation of $I_{\ell}^{0}(j)$ to
\begin{align*}
I_{\ell}^{0}(j) &= N_{\ell,0}^2 (\varphi_2 - \varphi_1) P_{\ell}(\cos\theta_j) 
\sum_{k=0}^{\ell} 
 {\scriptstyle \binom{\ell}{k} \binom{\ell + k}{k} }
\frac{(-1)^k}{2^k} \cdot \frac{1}{k+1}\times
\\&\hspace{5cm}
\left( (1 - \cos\vartheta_2)^{k+1} - (1 - \cos\vartheta_1)^{k+1} \right),    
\end{align*}
and confirms the corresponding formula in our statement.

Let $\ell\ge1$ and $1\le m\le \ell$ and consider the integral
\begin{align}\label{Imlj:first}
I_{\ell}^{m}(j)&=\frac{1}{m}N_{\ell,m}^2P^m_{\ell}(\cos\theta_j)\Big(\sin\big(m(\varphi_2-\phi_j)\big)-\sin\big(m(\varphi_1-\phi_j)\big)\Big)\times
\nonumber
\\
&\hspace{5cm}\times\int^{\cos\vartheta_1}_{\cos\vartheta_2}P_{\ell}^m(u)du.
\end{align}

All the attention is finally on the integral
\begin{equation*}
\int^{\cos\vartheta_1}_{\cos\vartheta_2}P_{\ell}^m(u)du.
\end{equation*}

By Lemma \ref{Lemma: integrals related to incomplete Beta function} and (\ref{Imlj:first}) we obtain
\begin{align*}
I^{m}_\ell (j)&=\frac{1}{m}N_{\ell,m}^2P^m_{\ell}(\cos\theta_j)\Big(\sin\big(m(\varphi_2-\phi_j)\big)-\sin\big(m(\varphi_1-\phi_j)\big)\Big)
\times
\\
&\hspace{0.6cm}\sum_{k=m}^{\ell}{\binom{\ell}{k}\binom{\ell+k}{k}}\frac{2(-1)^k k!}{(k-m)!}\cB_{k,m}(\cos\vartheta_2,\cos\vartheta_1),
\end{align*}
as claimed, and the proof is complete.
\end{proof}

\end{document}